# Spontaneous Fruit Fly Optimisation for truss weight minimisation: Performance evaluation based on the "no free lunch" theorem

Uche Onyekpe*, Stratis Kanarachos, Michael E. Fitzpatrick

onyekpeu@uni.coventry.ac.uk*, stratis.kanarachos@coventry.ac.uk, ab6856@coventry.ac.uk,

*Faculty of Engineering, Environment and Computing, Coventry University, Priory Street, Coventry*

*CV1 5FB, United Kingdom*

**Abstract**

*In recent years, researchers have presented various optimisation algorithms for truss design. The "no free lunch" theorem implies that no optimisation algorithm fits all problems; therefore, the interest is not only in their accuracy and convergence rate but also the tuning effort and population size to achieve the optimal result. The latter is particularly crucial for computationally-intensive or high-dimensional problems. Contrast-based Fruit-fly Optimisation Algorithm (c-FOA) proposed by Kanarachos et al. in 2017 is based on the efficiency of fruit flies in food foraging by olfaction and visual contrast. The proposed Spontaneous Fruit-fly Optimisation (s-FOA) enhances c-FOA and addresses the difficulty in solving nonlinear optimisation problems by presenting standard parameters and lean population size for use on optimisation problems. The s-FOA's performance is assessed using six benchmark problems. Comparison of the results obtained from documented literature and other investigated techniques demonstrates the competence and robustness of s-FOA in truss optimisation.*

## 1. Introduction

Trusses have found significant applications in modern engineering. Such applications range from use in transmission towers, offshore wind turbine supports, offshore oil and gas platforms, to microstructural applications such as the lattice structures of additive manufacturing [1, 2]. Truss optimisation aims to improve the performance of trusses while minimising the material resource [3]. The objective of the optimisation can be interpreted as a weight minimisation one, bounded by well-defined constraints. The constraints are the allowable stresses and displacements, as subject to high stress the truss members could fail through buckling or tension. There are many forms of optimisation, each with their unique design variables: this study focuses only on size optimisation. The design variable that is the most commonly investigated is the cross-sectional area of the truss member [4].

Optimisation algorithms are used in searching for the optimum solution to a problem. The application of optimisation algorithms to structures has proliferated in the last decade [5]. Many researchers have published on the applications of improved algorithms to truss weight minimisation problems. Kaveh and Mahdavi proposed a Multi-Objective Colliding Bodies Optimisation (MOCBO) algorithm for the optimisation of trusses bounded by an allowable stress limit [6]. A genetic programming methodology was used by Assimi *et al.* in the optimisation of the size and topology of trusses [7]. Another approach was taken by Cheng *et al.*, proposing a Hybrid Harmony Search (HHS) algorithm for the design of truss structures with stress limits [8]. Tejani *et al.* made use of the Improved Passing Vehicle Search (IPVS), Improved Heat Transfer Search (IHTS), Improved Water Wave Optimization (IWWO) and Improved Heat Transfer Search (IHTS) to optimise the topology of truss structures with displacement, stress and kinematic stability constraints [9]. An adaptive Symbiotic Organism Search (SOS) was utilised by Tejani *et al.* in truss structural optimisation with frequency constraints [10]. A development of PSO was presented by Kaveh and Zolghadr: Democratic PSO (Particle Swarm Optimisation) algorithm for the optimisation of truss layout and size with frequency constraints [11]. Multi-Class Teaching-Learning-Based Optimisation algorithm (MC-TLBO) was utilised by Farshchin *et al.* for truss design with frequency constraints [12]. Rajan used a Genetic Algorithm (GA) to optimise the shape, size and topology of truss structures [13].

Through all these research studies, the efficiency of optimisation algorithms in solving structural design problems has been established. However, according to the "no free lunch theorem", there exists no single algorithm to solve all optimisation problems. Hence the need to research lean algorithms [10].

In 2011, Pan proposed the FOA algorithm, a population-based technique that mimics the foraging activities of fruit flies [14]. Fruit-flies, compared to other species, possess a better sense of smell and vision which they use to find food efficiently. The algorithm has a framework which is simple, easy to understand, and is easily implementable in tackling various optimisation problems [15]. However, it is characterised by premature convergence (reduced accuracy) and is easily trapped in local optima.

FOA has been applied successfully to a variety of problems. In 2011, Pan applied the FOA algorithm to optimise the General Regression Neural Network and Multiple Regression utilised in modelling the financial distress problem of Taiwan's enterprise [14]. Lu *et al*. in 2015 proposed an adaptive fruit fly optimisation algorithm based on velocity variable (VFOA). The algorithm utilised the particle velocity concept from PSO on FOA to improve its convergence speed and accuracy. The improved algorithm was used to solve 13 mathematical benchmark problems [16]. As another improvement, Kanarachos *et al*. modified the FOA algorithm by including a visual contrast phase. The modification was based on biological discoveries on the complexity of the fruit foraging activities of fruit flies, thus improving its exploration capabilities. The modified algorithm was applied for the first time to solve truss design problems with stress, displacement or frequency constraints [17]. The algorithm was also used to improve the shock performance of vehicle suspension systems to mitigate damages caused by potholes in the UK [18]. Since then it was applied successfully in a range of problems. Wu *et al*. solved 33 mathematical benchmark functions by modifying the FOA to improve its exploration capabilities. A normal cloud generator was introduced to generate new positions of the swarm based on parameters such as possible food position, search range and search stability. It was inspired by the fact that fruit flies are characterised by fuzziness and randomness as they fly towards the food source. A cloud model is a tool used to synthesise the randomness and fuzziness in the algorithm [19]. Mitic *et al*. in 2015 presented the chaotic fruit fly optimisation algorithm to improve the explorative strategy of the algorithm. It does so by using the theory of chaos to relocate the fruit flies. The improved algorithm was used to solve ten one-dimensional benchmark mathematical problems [20]. With the aim of diversifying the solutions to avoid local optima or premature convergence, Yuan *et al*. introduced a Multi-Swarm Fruit Fly Optimisation Algorithm (MSFOA). The enhanced algorithm was used to identify parameters of a synchronous generator and solve six non-linear mathematical functions. In MSFOA, the swarm is divided into several sub-swarms, and the sub-swarms independently explore the search space to find the global optima [21].

FOA solves optimisation problems in two basic phases: the osphresis phase and the vision phase. The fruit fly makes use of its olfactory organ to detect the odour of the food source and then uses its vision capabilities to fly towards the food direction. Metaheuristic optimisation algorithms are characterised by two vital properties: exploration and exploitation. Exploration makes sure the algorithm visits a broader region of the search space (non-visited) for promising solutions. Exploitation aims to search intensively already visited regions of the search space for better solutions. Ensuring a good balance between exploration and exploitation is imperative to improve the performance of an optimisation algorithm and thus defines its success or failure. An unsuitable balance could lead to premature convergence, local optima entrapment and possibly stagnation [22]. Different problems require a different balance between exploitation and exploration, thus the need to adapt the parameters of the algorithm. In many cases, larger population sizes may compensate this problem. However, the requirement for large population sizes becomes problematic in higher dimensional problems.

The performance of optimisation algorithms significantly relies on certain unique parameters. Although many have tried, tuning of the algorithms to achieve the optimal result is not convenient. Algorithms such as GA, DE and PSO have their performances dependent on a number of parameters not known beforehand. The need for algorithms with standard set of parameters to achieve optimum results has therefore become necessary.

This study aims to present the Spontaneous Fruit fly Optimisation Algorithm (s-FOA) as an optimisation technique with a good exploration-exploitation balance characterised by fewer and standard tuning parameters. To this end, six benchmark truss design problems were solved with the s-FOA. The results were compared to several state of the art optimisation algorithms to establish the robustness of the algorithm.

## 2. Spontaneous Fruit Fly Optimisation Algorithm

Recently, a study by Van Breugel & Dickinson [23] showed that fruit flies exhibit a more complex search behaviour than the one modelled by Pan. Contrast-based Fruit Fly Optimisation (c-FOA) modelled for the first time the food search process of the Fruit Fly. The proposed Spontaneous Fruit Fly optimisation (s-FOA) enhances c-FOA and is described in detail in the following section. Its flowchart is provided in Figure 1.

2.1 Swarm generation, selection and termination

The algorithm starts by arbitrarily defining the position $(X_0, Y_0)$ of the first fruit fly in a coordinate system. Additional $N{-}1$ fruit flies are located, randomly, in the vicinity of $(X_0, Y_0)$ according to Eq. (1).

$$X_{ij}[k] = X_{0j}[k] \cdot (1 + M \cdot (2 \cdot rand_{N_{res}} - 1), j{=}1,2,...,m \text{ and } i{=}1,...,N$$

$$Y_{ij}[k] = Y_{0j}[k] \cdot (1 + M \cdot (2 \cdot rand_{N_{res}} - 1), j{=}1,2,...,m \text{ and } i{=}1,...,N$$

(1)

Where $k = 1,2,...,K_{max}$ is the iteration number, $m$ is the number of optimisation variables, $N$ is the size of the swarm and $rand_{N_{res}}$ is a random number, sampled from a uniform discrete distribution defined in the interval $[1, N_{res}]$. $M$ is a scaling parameter that defines how coarse or fine the search strategy is.

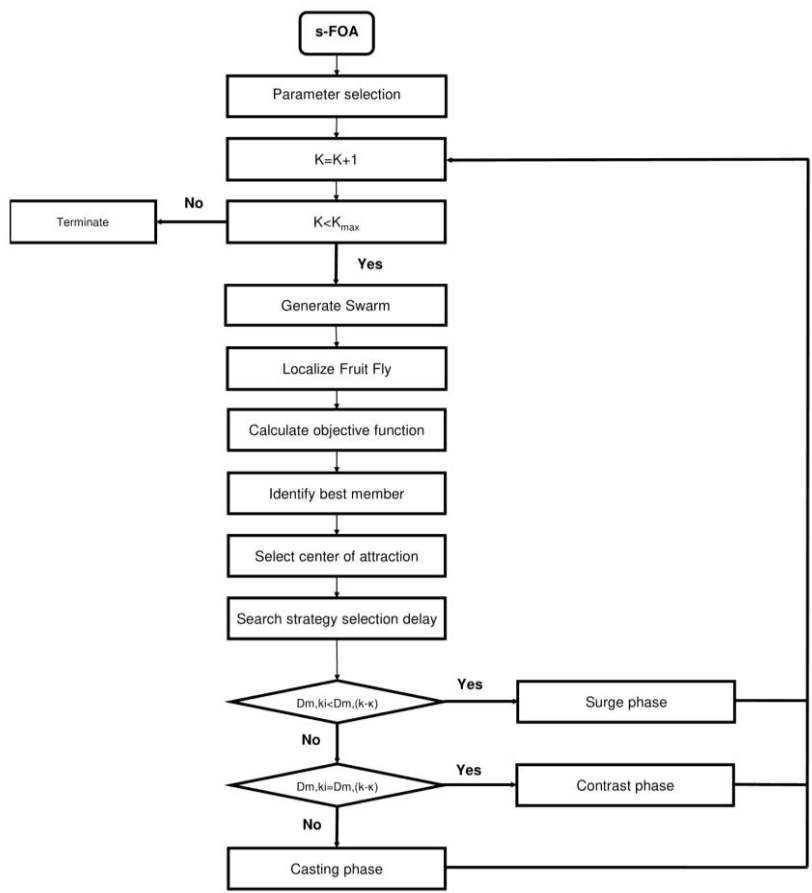

*Figure 1. Flowchart of Spontaneous Fruit Fly Optimisation Algorithm (s-FOA).*

Each fruit fly is assigned values $DI_{ij}$ based on how close each fruit fly parameter $(X_{ij}[k], Y_{ij}[k])$ is to the origin of the coordinate system:

$$D_{ij}[k] = \sqrt{X_{ij}^2[k] + Y_{ij}^2[k]} \tag{2}$$

$$DI_{ij}[k] = \frac{1}{D_{ij}[k]} \tag{3}$$

For each fruit fly at $\mathbf{d_i}[k]$ an objective function value $Dm_i[k]$ is assigned, $Dm_i[k] = f(\mathbf{d_i}[k])$.

In the following, the fruit flies are ranked based on their objective function values, and the fruit fly $\mathbf{d}^*[k]$ that achieves the lowest objective function value $Dm_i^*[k]$ at position $(X_i^*[k], Y_i^*[k])$ is identified. In case the objective function value $Dm_i^*[k]$ is lower than the previous centre of attraction $D_0[k]$, then $Dm_i^*[k]$ becomes the new centre of attraction $\mathbf{d}_0[k]$ $(X_0[k], Y_0[k])$.

$$if \ Dm_i^* < D_{m,k0}$$

$$then \ X_0[k] = X_i^*[k] \ and \ Y_0[k] = Y_i^*[k] \tag{4}$$

The algorithm terminates when the maximum number $K_{max}$ of iterations is reached.

In [23], it was shown that fruit flies exhibit three distinct behaviours when searching for food. First, when fruit flies find a better food source, they surge to that location. Second, fruit flies do not change immediately their search strategy when an improvement in food location is not achieved. Instead, they persist for a limited, constant amount of time. Third, fruit flies are attracted by visually contrasting objects. This means that they search also for food by using their eyes, not just smell. The three behaviours were included in s-FOA.

The delay in changing the search behaviour is represented by κ. The fruit flies can change their search strategy only every κ iterations. If the best objective function value $\mathbf{d}^*[k]$ improves over the last $\kappa$ iterations the swarm enters the "surging" phase, during which the fruit flies surge towards the attraction point $\mathbf{d}_0[k]$:

$$if \ (Dm_i^*[k] < Dm_i^*[k - \kappa])$$

$$M[k + 1] = c \cdot M[k] \tag{5}$$

In case the best objective function value does not change over the last $\kappa$ iterations the swarm enters the "visual contrast" phase, during which the fruit flies are attracted by an arbitrarily selected point $\mathbf{D}_{i,rand}[k]$:

$$if \ (Dm_i[k] = Dm_0[k - \kappa]$$

$$X_0[k] = X_{i,rand}[k] \ and \ Y_0[k] = Y_{i,rand}[k]] \tag{6}$$

where $k$ is the current iteration.

When a fruit fly does not improve its performance, then the swarm enters the "casting" phase. This behaviour is modelled based on [24]. There it was shown that fruit flies have memory that allows them to make decisions based on how good or bad the memory was.

$$if \ Dm_i[k] > Dm_i[k - 1]$$

$$then \ X_i[k] = X_i[k - 1] \ and \ Y_i[k] = Y_i[k - 1] \tag{7}$$

2.2 Centre of attraction and spontaneous positioning

In most FOA versions, the centre of attraction, the position towards which the fruit fly swarm moves, is usually the fruit fly that achieves the best performance. Although this is a natural step to take, the practice has shown that other algorithms, where the population moves towards random positions, perform much better in some truss optimisation problems. On the downside, this comes at the cost of significant performance degradation (slow convergence) for other problems.

A recent biological study showed that fruit flies have "free will" and can choose where they fly to, in the absence of a stimulus. This finding has triggered the development of the proposed Spontaneous Fruit Fly Optimisation (s-FOA)[1]. In s-FOA the fruit flies are attracted by a centre of attraction, which is defined by the combination of the positions of the fruit fly with the best performance and an arbitrarily selected fruit fly from the swarm.

Specifically, the positions of the fruit fly swarm in the next iterations are calculated according to:

$$X_{ij}[k] = X_i^*[k] + 0.5 \cdot \left(X_0[k] + X_{i,rand}[k]\right) \cdot M \cdot (2 \cdot rand_{N_{res}} - 1)$$
$$Y_{ij}[k] = Y_i^*[k] + 0.5 \cdot \left(Y_0[k] + Y_{i,rand}[k]\right) \cdot M \cdot (2 \cdot rand_{N_{res}} - 1)$$

(8)

where $j$=1,2,...,$m$ and $i$=1,...,$N$ and $(X_{i,rand}[k], Y_{i,rand}[k])$ is a randomly selected point for each $(X_{ij}[k], Y_{ij}[k])$.

The following pseudocode summarises s-FOA.

| | |
|---|---|
| 1 | **begin** |
| 2 | Select initial design vector $\mathbf{d}_0$=[DI$_{01}$, DI$_{02}$, ..., DI$_{0m}$], m is the number of design variables |
| 3 | Generate initial fruit fly swarm $\mathbf{d}_i$, i=1,2,..., N in the vicinity of $\mathbf{d}_0$ using uniform discrete distribution [1, N$_{res}$] |
| 4 | Calculate the objective function value Dm$_i$ at $\mathbf{d}_i$, Dm$_i$ = f($\mathbf{d}_i$) |
| 5 | Rank the fruit flies and find the one with the best performance |
| 6 | Dm$^*$ = f($\mathbf{d}^*$) = min(Dm$_i$) |
| 7 | If Dm$^*$ < Dm$_0$ then $\mathbf{d}_0$= $\mathbf{d}_i^*$ |
| 8 | **while (k< K)** |
| 9 | Increment k |
| 10 | Reposition the fruit fly swarm $\mathbf{d}_i$[k], near $\mathbf{d}_0$[k] using uniform discrete distribution [1, N$_{res}$] |
| 11 | Calculate the objective function value Dm$_i$[k] = f($\mathbf{d}_i$[k]) |
| 12 | Rank the fruit flies and find the best: |
| 13 | Dm$^*$[k] = f($\mathbf{d}^*$[k]) = min(Dm$_i$[k]) |
| 14 | If Dm$^*$[k] < Dm$_0$[k] then $\mathbf{d}_0$[k+1]=$\mathbf{d}^*$[k] |
| 15 | Increment response time t[k] = t[k−1]+1 |
| 16 | **If** (t[k]>delay time) |
| 17 | **If** (Dm$^*$[k] < Dm$_0$[k − κ]) |
| 18 | reduce the search radius M[k + 1] = c · M[k] |
| 19 | (surging phase) |
| 20 | **else if** (Dm$^*$[k] = Dm$_0$[k − κ]) |
| 21 | a random candidate, $\mathbf{d}_{rand}$[k] , |

[1] https://io9.gizmodo.com/the-crazy-device-that-shows-fruit-flies-have-free-will-1459261376

| | |
|---|---|
| 22 | $D_{rand}[k] = f(\mathbf{d_{rand}}[k])$, becomes the new attraction point<br>$\mathbf{d}_0[k+1] = \mathbf{d_{rand}}[k]$ |
| 23 | (contrast-based vision phase) |
| 24 | **end if** |
| 25 | Initialise response time t[k]=0 |
| 26 | **end if** |
| 27 | **End while** |
| 28 | Post process results and visualisation |
| 29 | **end** |

Figure 2 Pseudo-code of s-FOA.

## 3. Mathematical formulation of the truss design problem

The goal of truss optimisation is to minimize the weight of the structure such that constraints on its performance are satisfied. In this study, the design variable is chosen as the cross-sectional area of the members of the truss structures. The mathematical formulation of the truss optimisation problem is as shown below:

To minimize the weight of the truss,

$$Minimise\ W = \sum_{i=1}^{m} A_i\, \rho_i L_i \ \ i = 1,2,3,\dots\dots,m \tag{9}$$

Where $W$ is the weight of the truss structure consisting of $m$ members and $A_i, \rho_i$ and $L_i$ are respectively the cross-sectional area, material density and length of the $i$th truss member. The design constraints are defined as follows:
Subject to:

$$\sigma_{min} \leq \sigma_i \leq \sigma_{max}, \qquad i = 1,2,3,\dots\dots,m$$

$$\delta_{min} \leq \delta_j \leq \delta_{max}, \ j = 1,2,3,\dots\dots,n \tag{10}$$

$$\sigma_k^b \leq \sigma_k \leq 0, \qquad k = 1,2,3,\dots\dots,mc$$

$$A_i \in S = \{A_1, A_2, A_3, \dots\dots, A_d\}$$

Where $\sigma_i$ and $\delta_j$ are the $i$th member allowable stress and $j$th nodal displacement respectively; $\sigma_k$ is the allowable buckling stress of the $i$th member under compression; $S$ is a set of discrete cross-sectional area; $m, mc$ and $n$ are the total number of members, members subject to compression and nodes in the truss structure respectively.

In this work, the penalty approach is adopted for the transformation of the constrained optimisation problem to an unconstrained problem. Consequently, the mathematical formulation of the truss optimisation problem becomes:

$$f_{penalty} W = W + \ \lambda \sum_{r=1}^{nin} [\max(0, g_r)]^2 + \lambda \sum_{s=1}^{neq} [\max(0, |h_s|) - \epsilon]^2 \tag{11}$$

Where $nin$ and $neq$ are the number of equality and inequality constraints; $g_r$ and $h_s$ are the $r$th equality and $s$th inequality constraints; $\lambda$ is the penalty value chosen as $10^5$ and is the small positive tolerance for the equality constraints chosen as $10^{-6}$.

## 4. Benchmark truss design problems and discussions

In this section, six benchmark truss problems consisting of 10, 15, 25, 52, 72 and 200 members are studied as a size optimisation problem constrained by nodal displacement and/or allowable stresses in the truss members. For each of the six truss optimisation problems, the performance of the s-FOA algorithm is investigated and compared to other implemented algorithms such as TLBO (Teaching Learning Based Optimisation), GA, DE (Differential Evolution) and PSO. A comparison is also made to the results reported in published literatures. All the benchmark problems were coded and the optimisation algorithms implemented in MatLab R2017a. A statistical stochastic performance analysis of each algorithm is also investigated by conducting 30 independent runs per problem. The DE, PSO, TLBO and GA algorithm used in this study are sourced from [25], [26], [27] and [28] respectively.

The operational parameters of the DE algorithm were selected according to [29]. The population size (N), cross-over probability (Cr), mutation factor (F), parameters chosen, and the maximum number of functional evaluations are as displayed on Table 1. The population members were created by a uniform random distribution. All variables were treated internally as floating variables by DE. The algorithm was terminated when the maximum number of functional evaluations is reached.

*Table 1 Selected parameters for DE algorithm*

| DE | | | | | | |
|---|---|---|---|---|---|---|
| Truss Problem | 10 Bar | 15 Bar | 25 Bar | 52 Bar | 72 Bar | 200 Bar |
| Number of Fitness Evaluations | 2000 | 2000 | 2000 | 14000 | 10500 | 48000 |
| Number of unknowns | 10 | 15 | 8 | 12 | 16 | 29 |
| Population size | 12 | 28 | 12 | 28 | 37 | 75 |
| Crossover probability (Cr) | 0.2368 | 0.9426 | 0.2368 | 0.9426 | 0.9455 | 0.8803 |
| Mutation factor (F) | 0.6702 | 0.6607 | 0.6702 | 0.6607 | 0.6497 | 0.4717 |
| Number of generations | 167 | 71.429 | 167 | 500 | 283.78 | 600 |

For the PSO algorithm, the number of agents in the population (N), inertia weight, maximum velocity chosen and maximum number of functional analysis as shown on Table 2 were obtained from [30]. The initial swarm members were created by a random distribution. Subsequently, new members of the swarm were created using:

$$v_{i+1} = \omega \cdot v_i + c_1 \cdot r_p \cdot (p_i - x_i) + c_2 \cdot r_g \cdot (p_g - x_i) \qquad (12)$$

$$x_{i+1} = x_i + v_{i+1}$$

Where parameters $c_1$ and $c_2$ were chosen to be 1.5 and 2.0 respectively. The inertia weight was selected within the range [0.8, 1.2] to prevent weak exploration and local optima entrapment. The algorithm was terminated when the number the maximum of functional evaluations is reached.

*Table 1 Selected parameters for PSO algorithm*

| PSO | | | | | | |
|---|---|---|---|---|---|---|
| Truss Problem | 10 Bar | 15 Bar | 25 Bar | 52 Bar | 72 Bar | 200 Bar |
| Number of Fitness Evaluations | 2000 | 2000 | 2000 | 14000 | 10500 | 48000 |
| Number of unknowns | 10 | 15 | 8 | 12 | 16 | 29 |
| Population size | 12 | 28 | 37 | 28 | 37 | 75 |
| Number of generations | 167 | 71 | 54 | 500 | 284 | 600 |
| Inertia weight | 0.8 | | | | | |
| Vmax | Upper bound Cross-Sectional Area | | | | | |

To ensure optimal performance of the GA algorithm, tuning parameters were selected from [31]. The size of the population (N), cross-over probability (Cr), mutation factor (F) as well as the maximum number of structural analyses selected are shown in Table 3. Initial population members were generated by a uniform random distribution. The parents in each generation were selected stochastically and then weighted by the crossover operator for the creation of new members. Mutation diversified the population members through random selection.

*Table 2 Selected parameters for GA algorithm*

| GA | | | | | | |
| --- | --- | --- | --- | --- | --- | --- |
| Truss Problem | 10 Bar | 15 Bar | 25 Bar | 52 Bar | 72 Bar | 200 Bar |
| Number of Fitness Evaluations | 2000 | 2000 | 2000 | 14000 | 10500 | 48000 |
| Number of unknowns | 10 | 15 | 8 | 12 | 16 | 29 |
| Population size | 37 | 37 | 37 | 28 | 37 | 75 |
| Number of generations | 54 | 54 | 54 | 500 | 284 | 600 |

The TLBO algorithm utilised is that from [27] with a maximum number of analysis as shown on Table 4. The class members were created randomly using a uniform distribution. The best student is selected as the class teacher in the teacher phase. In the learner's phase, the class members improve the individual and class performance through student-student interaction. Successive implementation of both phases goes on until the number of maximum of functional evaluations is reached.

*Table 4 Selected parameters for TLBO algorithm*

| TLBO | | | | | | |
| --- | --- | --- | --- | --- | --- | --- |
| Truss Problem | 10 Bar | 15 Bar | 25 Bar | 52 Bar | 72 Bar | 200 Bar |
| Number of Fitness Evaluations | 2000 | 2000 | 2000 | 14000 | 10500 | 48000 |
| Number of unknowns | 10 | 15 | 8 | 12 | 16 | 29 |
| Population size | 28 | 28 | 12 | 28 | 37 | 75 |
| Number of generations | 36 | 36 | 84 | 250 | 142 | 300 |

The parameters for the c-FOA algorithm as shown on Table 5 are selected as: Population size=50, $\kappa$=320, $M$=0.95, $N_{res}$=50 and $c$=0.92.

*Table 5 Selected parameters for cFOA algorithm*

| cFOA | | | | | | |
| --- | --- | --- | --- | --- | --- | --- |
| Truss Problem | 10 Bar | 15 Bar | 25 Bar | 52 Bar | 72 Bar | 200 Bar |
| Number of Fitness Evaluations | 2000 | 2000 | 2000 | 14000 | 9600 | 45000 |
| Number of unknowns | 10 | 15 | 8 | 12 | 16 | 29 |
| Population size | 10 | 15 | 8 | 20 | 16 | 60 |
| Number of generations | 200 | 130 | 250 | 700 | 600 | 800 |

The parameters for the s-FOA algorithm as shown on Table 5 are selected as: Population size=50, $\kappa$=5, $M$=0.95, $N_{res}$=10 and $c$=0.9.

*Table 6 Selected parameters for s-FOA algorithm*

| s-FOA | | | | | | |
| --- | --- | --- | --- | --- | --- | --- |
| Truss Problem | 10 Bar | 15 Bar | 25 Bar | 52 Bar | 72 Bar | 200 Bar |
| Number of Fitness Evaluations | 2000 | 1950 | 2000 | 14000 | 9600 | 45000 |
| Number of unknowns | 10 | 15 | 8 | 12 | 16 | 29 |
| Population size | 10 | 15 | 8 | 20 | 16 | 60 |
| Number of generations | 200 | 130 | 250 | 700 | 600 | 800 |

## 4.1 : Benchmark 1: 10-Bar planar truss

### 4.1.1 Benchmark 1: Problem Description

The first truss problem is a 10-bar planar structure as shown in Figure 3. The design was examined by Sadollah *et al.* [33], Li *et al.* [34], Ho-Huu *et al.* [35], Do and Lee [36] and many others as a size optimisation problem with ten design variables. The design variable is defined as the cross-sectional area of the truss members and is chosen from a discrete set of data S= {1.62, 1.80, 1.99, 2.13, 2.38, 2.62, 2.88, 2.93, 3.09, 3.13,3.38, 3.47, 3.55, 3.63, 3.84, 3.87, 3.88, 4.18, 4.22, 4.49, 4.59, 4.80,4.97, 5.12, 5.74, 7.22, 7.97, 11.50, 13.50, 13.90, 14.20, 15.50, 16.00,16.90, 18.80, 19.90, 22.00, 22.90, 26.50, 30.00, 33.50}in². Each member of the truss has a material density and a modulus of elasticity as defined in Table 6. Two vertical loads of 10000 lbs each are vertically applied to the truss at nodes 2 and 4 with each member and node constrained as also shown in Table 6.

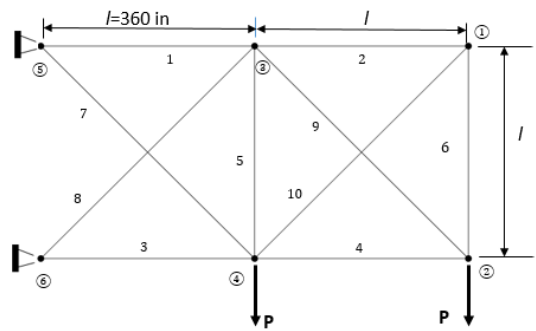

*Figure 3 Benchmark Case 1: The 10-bar planar truss structure*

*Table 7 Material properties and design constraints of the benchmark truss problems*

| Truss Problem | 10 Bar | 15 Bar | 25 Bar | 52 Bar | 72 Bar | 200 Bar |
|---|---|---|---|---|---|---|
| Material Density ρ | 0.1 lb/in3 | 7800 kg/m3 | 0.1 lb/in3 | 7800 kg/m3 | 0.1 lb/in3 | 0.283 lb/in3 |
| Modulus of Elasticity E | 104 ksi | 200 GPa | 104 ksi | 2.07 x105 MPa | 104 ksi | 30000 ksi |
| Nodal Displacement Constraint δ | ± 2 inches the x and y direction | ± 10 mm the x and y direction | ± 0.35 inches the x and y direction | N/A | ± 0.25 inches the x and y direction | N/A |
| Stress Constraint σ | ±25ksi in tension and compression | ±120MPa in tension and compression | ±40ksi in tension and compression | ±180MPa in tension and compression | ±25ksi in tension and compression | ±10 ksi in tension and compression |

### 4.1.2 Benchmark 1: Results and discussion

The optimal results obtained from the run of all five algorithms and the MBA [33], HPSO [34], aeDE [35], SOS [36] and mSOS [36] are presented in Tables 7 and 8. The results highlight the size variables, best weight, mean weight, standard deviation and number of functional evaluations for all algorithms in this study and in reported literature. From the results, it can be observed that the s-FOA produces the best continuous solution without any constraint violation as 5421.2 lb compared to the other algorithms implemented in the investigation. Relatively, the s-FOA also proves to be the most consistent of the implemented algorithms by recording the lowest standard deviation of 41.62 lb. It further displays the best mean weight of 5506.2 lb.

The best discrete optimum result is achieved by s-FOA,DE, PSO, TLBO, SOS, mSOS, and aeDE with a weight of 5490.738 lb. However, the maximum number of functional analyses used in this research is small compared to that implemented in the reported literatures.

*Table 8 Benchmark Case 1: Performance rankings of the algorithms in solving the 10 bar truss problem*

|  | PERFORMANCE RANKINGS | | | | | |
|---|---|---|---|---|---|---|
|  | DE | PSO | GA | TLBO | cFOA | s-FOA |
| Best Value | 4 | 3 | 6 | 5 | 2 | 1 |
| Mean Value | 5 | 2 | 6 | 4 | 3 | 1 |
| Number of Functional Evaluations per 10,000 | 0.2 | 0.2 | 0.2 | 0.2 | 0.2 | 0.2 |
| Standard deviation | 5 | 3 | 6 | 4 | 2 | 1 |
| Number of tuning parameters changed from one benchmark problem to another | 5 | 3 | 3 | 3 | 3 | 3 |
| Total Score | 16.2 | 9.2 | 18.2 | 13.2 | 10.2 | 6.2 |

*Table 9 Optimal results of the 10-Bar planar truss structure optimisation obtained from different algorithms from previous studies*

| Area group | PREVIOUS STUDY | | | | | |
|---|---|---|---|---|---|---|
|  | HPSO [34] | MBA [33] | aeDE [35] | DE [36] | SOS [36] | mSOS [36] |
| in$^2$ | Discrete optimum | Discrete optimum | Discrete optimum | Discrete optimum | Discrete optimum | Discrete optimum |
| A1 | 30 | 30 | 33.5 | 33.5 | 33.5 | 33.5 |
| A2 | 1.62 | 1.62 | 1.62 | 1.62 | 1.62 | 1.62 |
| A3 | 22.9 | 22.9 | 22.9 | 22.9 | 22.9 | 22.9 |
| A4 | 13.5 | 16.9 | 14.2 | 14.2 | 14.2 | 14.2 |
| A5 | 1.62 | 1.62 | 1.62 | 1.62 | 1.62 | 1.62 |
| A6 | 1.62 | 1.62 | 1.62 | 1.62 | 1.62 | 1.62 |
| A7 | 7.97 | 7.97 | 7.97 | 7.97 | 7.97 | 7.97 |
| A8 | 26.5 | 22.9 | 22.9 | 22.9 | 22.9 | 22.9 |
| A9 | 22 | 22.9 | 22 | 22 | 22 | 22 |
| A10 | 1.8 | 1.62 | 1.62 | 1.62 | 1.62 | 1.62 |
| Weight (lb) | 5531.98 | 5507.75 | 5490.738 | 5490.738 | 5490.738 | 5490.738 |
| Mean weight (lb) | - | - | - | - | - | - |
| Standard Deviation | - | - | - | - | - | - |
| Number of Functional Evaluations | 50000 | 3600 | 2380 | 300,000 | 300,000 | 300,000 |

*Table 10 Optimal results of the 10-Bar planar truss structure optimisation obtained from different algorithms in this study*

| Area group | THIS STUDY | | | | | | | | | | | |
| | TLBO | | PSO | | DE | | GA | | cFOA | | s-FOA | |
| in² | Continuous Optimum | Discrete optimum | Continuous Optimum | Discrete optimum | Continuous Optimum | Discrete optimum | Continuous Optimum | Discrete optimum | Continuous Optimum | Discrete Optimum | continuous | Discrete optimum |
| A1 | 32.7009 | 33.5 | 31.7749 | 33.5 | 33.5 | 33.5 | 32.7403 | 30 | 32.554386 | 33.5 | 31.6141 | 33.5 |
| A2 | 1.62 | 1.62 | 1.6295 | 1.62 | 1.6462 | 1.62 | 3.0609 | 2.93 | 1.613358 | 1.62 | 1.6218 | 1.62 |
| A3 | 23.8969 | 22.9 | 23.906 | 22 | 22.4672 | 22.9 | 25.5496 | 22.9 | 22.782222 | 22.9 | 23.3547 | 22.9 |
| A4 | 15.0507 | 14.2 | 15.146 | 14.2 | 15.4151 | 14.2 | 19.8758 | 18.8 | 15.63543 | 14.2 | 14.8084 | 14.2 |
| A5 | 1.636 | 1.62 | 1.6201 | 1.62 | 1.6209 | 1.62 | 1.8382 | 1.8 | 1.615788 | 1.62 | 1.1627 | 1.62 |
| A6 | 1.6396 | 1.62 | 1.62 | 1.62 | 1.6203 | 1.62 | 2.4316 | 2.38 | 1.622592 | 1.62 | 1.2364 | 1.62 |
| A7 | 8.146 | 7.97 | 8.6446 | 7.97 | 7.9847 | 7.97 | 11.7973 | 7.97 | 8.411202 | 7.97 | 8.269 | 7.97 |
| A8 | 21.7906 | 22.9 | 22.8267 | 22.9 | 23.3522 | 22.9 | 19.1187 | 18.8 | 22.61358 | 22.9 | 22.8569 | 22.9 |
| A9 | 22.2357 | 22 | 21.1247 | 22 | 20.9034 | 22 | 17.9177 | 16.9 | 21.460788 | 22 | 21.534 | 22 |
| A10 | 1.6508 | 1.62 | 1.62 | 1.62 | 1.62 | 1.62 | 2.6401 | 2.38 | 1.618056 | 1.62 | 1.6257 | 1.62 |
| Weight (lb) | 5495.8 | 5490.738 | 5485.2 | 5490.738 | 5487.8 | 5490.738 | 5698.5 | 5526.962 | 5484.2 | 5491.717 | 5421.2 | 5490.738 |
| Mean weight (lb) | 5570.5 | | 5539.8 | | 5688.9 | | 6223.7 | | 5555 | | 5506.2 | |
| Standard Deviation | 130.9274 | | 61.3838 | | 253.0152 | | 359.0314 | | 46.1695 | | 41.62 | |
| Number of Functional Evaluations | 2000 | | 2000 | | 2000 | | 2000 | | | | 2000 | |

4.2 : Benchmark 2: 15-Bar planar truss

4.2.1 Benchmark 2: Problem description

Another benchmark truss problem considered is the 15-bar planar truss. The problem has been investigated by researchers including Li *et al.* [34]. Each member of this truss problem has its material property and design constraints defined on Table 6. The problem has 15 design variables chosen from a discrete set of data S = {113.2, 143.2, 145.9, 174.9, 185.9, 235.9, 265.9, 297.1, 308.6, 334.3, 338.2, 497.8, 507.6, 736.7, 791.2, 1063.7} mm$^2$. The load cases acting on the structure are: Case 1: P1 = 35 kN, P2 = 35 kN, P3 = 35 kN; Case 2: P1 = 35 kN, P2 = 0 kN, P3 = 35 kN; Case 3: P1 = 35 kN, P2 = 35 kN, P3 = 0 kN Figure 4 shows the geometry and loading arrangement of the truss.

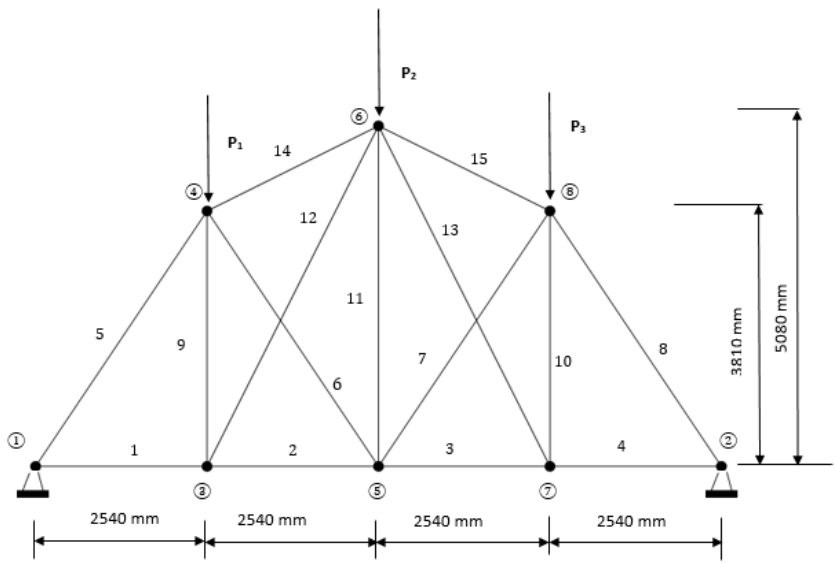

*Figure 4 Benchmark Case 2: The 15-Bar planar truss structure*

4.2.2 Benchmark 2: Results and discussion

The optimal design of the truss is obtained by taking all three load cases into consideration. A summary of the result is presented in Tables 9 and 10. From the analysis of the continuous results investigated in this study, the s-FOA produces the lightest truss weight of 86.57 kg compared to the other algorithms. s-FOA also records the lowest standard deviation of 1.70 kg and mean weight of 89.75 kg, yet again proving the stability and robustness of the algorithm. From the discrete solution analysis, the s-FOA, TLBO, DE, and HPSO [34] provide the best weight of 105.74 kg each compared to that of cFOA, GA and PSOPC [34] with weights of 107.35kg, 117.96 kg and 108.96 kg.

*Table 11 Benchmark Case 2: Performance rankings of the algorithms for the 15 bar truss problem*

|  | PERFORMANCE RANKINGS | | | | | |
|---|---|---|---|---|---|---|
|  | DE | PSO | GA | TLBO | cFOA | s-FOA |
| Best Value | 5 | 2 | 6 | 3 | 4 | 1 |
| Mean Value | 5 | 2 | 6 | 3 | 4 | 1 |
| Number of Functional Evaluations per 10,000 | 0.2 | 0.2 | 0.2 | 0.2 | 0.2 | 0.195 |
| Standard deviation | 5 | 4 | 6 | 3 | 1 | 2 |
| Number of tuning parameters changed from one benchmark problem to another | 5 | 3 | 3 | 3 | 3 | 3 |
| Total Score | 17.2 | 10.2 | 18.2 | 11.2 | 12.2 | 6.195 |

*Table 12 Optimal results of the 15-bar planar truss structure optimisation obtained from different algorithms*

| | | | | 15 Bar | | | | | | | | | | | |
|---|---|---|---|---|---|---|---|---|---|---|---|---|---|---|---|
| Area group | PREVIOUS STUDY | | | THIS STUDY | | | | | | | | | | | |
| | PSO [34] | PSOPC [34] | HPSO [34] | TLBO | | PSO | | DE | | GA | | cFOA | | s-FOA | |
| mm² | Discrete optimum | Discrete optimum | Discrete optimum | Continuous Optimum | Discrete optimum | Continuous Optimum | Discrete optimum | Continuous Optimum | Discrete optimum | Continuous Optimum | Discrete optimum | Continuous Optimum | Discrete Optimum | Continuous Optimum | Discrete optimum |
| A1 | 185.9 | 113.2 | 113.2 | 113.2 | 113.2 | 113.2 | 113.2 | 113.2 | 113.2 | 135.8017 | 113.2 | 122.58428 | 113.2 | 139.0209 | 113.2 |
| A2 | 113.2 | 113.2 | 113.2 | 120.7607 | 113.2 | 113.2 | 113.2 | 113.2 | 113.2 | 128.8207 | 113.2 | 166.90208 | 145.9 | 114.3433 | 113.2 |
| A3 | 143.2 | 113.2 | 113.2 | 113.2 | 113.2 | 113.2105 | 113.2 | 141.3093 | 113.2 | 113.2 | 113.2 | 123.8974 | 113.2 | 129.614 | 113.2 |
| A4 | 113.2 | 113.2 | 113.2 | 143.9061 | 143.2 | 113.2 | 113.2 | 113.2 | 113.2 | 176.152 | 174.9 | 134.06276 | 113.2 | 92.4504 | 113.2 |
| A5 | 736.7 | 736.7 | 736.7 | 526.5156 | 736.7 | 525.8146 | 736.7 | 539.219 | 736.7 | 673.1077 | 736.7 | 530.1722 | 736.7 | 528.5308 | 736.7 |
| A6 | 143.2 | 113.2 | 113.2 | 118.266 | 113.2 | 113.2 | 113.2 | 129.2336 | 113.2 | 123.0755 | 113.2 | 117.93176 | 113.2 | 81.6285 | 113.2 |
| A7 | 113.2 | 113.2 | 113.2 | 118.606 | 113.2 | 113.2 | 113.2 | 141.3216 | 113.2 | 113.2 | 113.2 | 121.53152 | 113.2 | 87.3564 | 113.2 |
| A8 | 736.7 | 736.7 | 736.7 | 531.2776 | 736.7 | 525.8101 | 736.7 | 527.013 | 736.7 | 788.847 | 736.7 | 526.06304 | 736.7 | 530.342 | 736.7 |
| A9 | 113.2 | 113.2 | 113.2 | 113.6073 | 113.2 | 113.2 | 113.2 | 124.7797 | 113.2 | 138.0283 | 113.2 | 152.6502 | 145.9 | 88.3752 | 113.2 |
| A10 | 114.2 | 113.2 | 113.2 | 120.7674 | 113.2 | 113.2 | 113.2 | 123.2264 | 113.2 | 169.4542 | 145.9 | 126.35384 | 113.2 | 110.3134 | 113.2 |
| A11 | 113.2 | 113.2 | 113.2 | 113.2 | 113.2 | 113.2 | 113.2 | 125.7362 | 113.2 | 202.6266 | 185.9 | 118.10156 | 113.2 | 83.0209 | 113.2 |
| A12 | 116.2 | 113.2 | 113.2 | 113.2 | 113.2 | 113.2114 | 113.2 | 113.2 | 113.2 | 117.936 | 113.2 | 112.9736 | 113.2 | 105.2534 | 113.2 |
| A13 | 117.2 | 185.9 | 114.2 | 113.2 | 113.2 | 113.2088 | 113.2 | 120.6281 | 113.2 | 196.3723 | 185.9 | 119.79956 | 113.2 | 103.5667 | 113.2 |
| A14 | 334.3 | 334.3 | 334.3 | 322.484 | 334.3 | 321.2347 | 334.3 | 528.2988 | 334.3 | 540.2041 | 507.6 | 338.64912 | 334.3 | 325.5745 | 334.3 |
| A15 | 334.3 | 334.3 | 334.3 | 328.466 | 334.3 | 321.2326 | 334.3 | 335.0375 | 334.3 | 476.7704 | 338.2 | 333.3174 | 334.3 | 344.8298 | 334.3 |
| Weight(kg) | 108.84 | 108.96 | 105.735 | 91.8692 | 106.329 | 90.093 | 105.735 | 99.1078 | 105.735 | 125.2476 | 117.956 | 95.2911 | 107.335 | 86.5694 | 105.735 |
| Mean weight (kg) | - | - | - | 94.1729 | | 91.9209 | | 105.5734 | | 152.0246 | | 96.544 | | 89.7546 | |
| Standard Deviation | - | - | - | 1.8892 | | 5.2483 | | 6.4822 | | 14.0832 | | 1.2195 | | 1.7082 | |
| Number of Functional Evaluations | 25000 | 25000 | 25000 | 2000 | | 2000 | | 2000 | | 2000 | | | | 1950 | |

## 4.3 : Benchmark 3 – 25-Bar planar truss

### 4.3.1 Benchmark 3: Problem description

The third benchmark truss illustrated in Figure 5 is the 25-bar space truss. The truss design problem previously investigated by Sadollah *et al.*, Li *et al.*, Wu *et al.* and Lee *et al.* [33, 34, 37, 38] has its material properties and design constraints defined in Table 6. The loading configuration of the truss is defined in Table 11 with each member categorised into a group of 8 representing the design variables and chosen from a discrete set of data S= {0.1, 0.2, 0.3, 0.4, 0.5, 0.6, 0.7, 0.8, 0.9, 1.0, 1.1, 1.2, 1.3, 1.4, 1.5, 1.6, 1.7, 1.8, 1.9, 2.0, 2.1, 2.2, 2.3, 2.4, 2.6, 2.8, 3.0, 3.2, 3.4} (in$^2$). The 8 design groups are as follows: (1) A1; (2) A2–A5; (3) A6–A9; (4) A10–A11; (5) A12–A13; (6) A14–A17; (7) A18–A21, and (8) A22–A25.

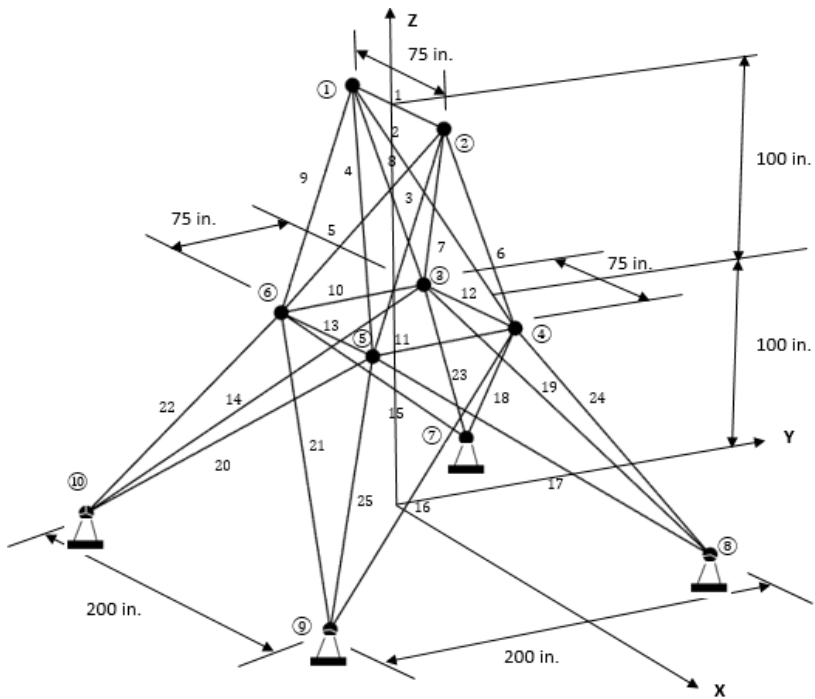

*Figure 5 Benchmark Case 3: The 25-Bar spatial truss structure*

*Table 13 Benchmark Case 3: Loading Configuration of the 25-bar truss*

| Nodes | Loads | | |
|---|---|---|---|
| | $P_x$ (kips) | $P_y$ (kips) | $P_z$ (kips) |
| 1 | 1 | −10 | −10 |
| 2 | 0 | −10 | −10 |
| 3 | 0.5 | 0 | 0 |
| 6 | 0.6 | 0 | 0 |

### 4.3.2 Benchmark 3: Results and discussion

A comparison of the optimal and statistical result obtained for the 25-bar truss design is as addressed in Table 12 and 13. The best continuous truss design weighing 483.9 lb is provided by s-FOA and cFOA. The s-FOA also gives the lowest standard deviation and second best mean weight of 9.97 lb and 492.9 lb respectively as compared to the PSO, DE, GA and TLBO algorithms. Therefore, it can be expressed that the s-FOA exhibits more robustness and accuracy over the compared algorithms from the view point of continuous optimum design.

From a discrete optimisation perspective, the design found by the s-FOA and cFOA proves to be the lightest truss weighing 483.87lb compared to the, MBA [33], HPSO [34], SOS [36], mSOS [36] and HS [38] of 484.85lb each, SGA [37] of 486.29 lb, DE of 485.57 lb, TLBO of 487.87 lb, PSO of 488.57 lb and GA of 489.49 lb. The performance of the s-FOA is further evidenced by the attainment of the solution under a significantly lower number of structural analysis of 1950 compared to 300,000 analyses required by the SOS [36] and mSOS [36], 40000 analyses by the SGA [37], 25000 analyses by the HPSO [34], 18734 analyses by the HS [38] and 3750 by the MBA [33].

*Table 14 Benchmark Case 3: Performance rankings of the algorithms for the 25-bar truss problem*

|  | PERFORMANCE RANKINGS | | | | | |
|---|---|---|---|---|---|---|
|  | DE | PSO | GA | TLBO | cFOA | s-FOA |
| Best Value | 5 | 4 | 6 | 3 | 1 | 1 |
| Mean Value | 5 | 4 | 6 | 1 | 3 | 2 |
| Number of Functional Evaluations per 10,000 | 0.2 | 0.2 | 0.2 | 0.2 | 0.2 | 0.2 |
| Standard deviation | 2 | 6 | 5 | 4 | 3 | 1 |
| Number of tuning parameters changed from one benchmark problem to another | 5 | 3 | 3 | 3 | 3 | 3 |
| Total Score | 15.2 | 14.2 | 17.2 | 9.2 | 10.2 | 7.2 |

*Table 15 Optimal results of the 25-bar space truss structure optimisation obtained from different algorithms from previous study*

| 25 Bar | | | | | | | |
|---|---|---|---|---|---|---|---|
| PREVIOUS STUDY | | | | | | | |
| **Area group** | SGA [37] | HS [38] | HPSO [34] | MBA [33] | DE [36] | SOS [36] | mSOS [36] |
| in² | Discrete optimum | Discrete optimum | Discrete optimum | Discrete optimum | Discrete optimum | Discrete optimum | Discrete optimum |
| A1 | 0.1 | 0.1 | 0.1 | 0.1 | 0.1 | 0.1 | 0.1 |
| A2 | 0.5 | 0.3 | 0.3 | 0.3 | 0.3 | 0.3 | 0.3 |
| A3 | 3.4 | 3.4 | 3.4 | 3.4 | 3.4 | 3.4 | 3.4 |
| A4 | 0.1 | 0.1 | 0.1 | 0.1 | 0.1 | 0.1 | 0.1 |
| A5 | 1.5 | 2.1 | 2.1 | 2.1 | 2.1 | 2.1 | 2.1 |
| A6 | 0.9 | 1 | 1 | 1 | 1 | 1 | 1 |
| A7 | 0.6 | 0.5 | 0.5 | 0.5 | 0.5 | 0.5 | 0.5 |
| A8 | 3.4 | 3.4 | 3.4 | 3.4 | 3.4 | 3.4 | 3.4 |
| Weight (lb) | 486.29 | 484.85 | 484.85 | 484.85 | 484.85 | 484.85 | 484.85 |
| Mean weight | - | - | - | - | - | - | - |
| Standard Deviation | - | - | - | - | - | - | - |
| Functional Evaluations | 40000 | 18734 | 25000 | 3750 | 300,000 | 300,000 | 300,000 |

*Table 16 Optimal results of the 25-bar space truss structure optimisation obtained from different algorithms in this study*

| | THIS STUDY | | | | | | | | | | | |
|---|---|---|---|---|---|---|---|---|---|---|---|---|
| **Area group** | **TLBO** | | **PSO** | | **DE** | | **GA** | | **cFOA** | | **s-FOA** | |
| in$^2$ | Continuous Optimum | Discrete optimum | Continuous Optimum | Discrete optimum | Continuous Optimum | Discrete optimum | Continuous Optimum | Discrete optimum | Continuous Optimum | Discrete Optimum | Continuous Optimum | Discrete optimum |
| A1 | 0.1012 | 0.1 | 0.1 | 0.1 | 0.1 | 0.1 | 0.1043 | 0.1 | 0.14331 | 0.1 | 0.14331 | 0.1 |
| A2 | 0.3621 | 0.4 | 0.3163 | 0.4 | 0.3864 | 0.4 | 0.7715 | 0.7 | 0.39318 | 0.4 | 0.39318 | 0.4 |
| A3 | 3.3967 | 3.4 | 3.4 | 3.4 | 3.4 | 3.4 | 3.345 | 3.4 | 3.4493 | 3.4 | 3.4493 | 3.4 |
| A4 | 0.1019 | 0.1 | 0.1 | 0.1 | 0.1 | 0.1 | 0.126 | 0.1 | 0.25793 | 0.2 | 0.25793 | 0.2 |
| A5 | 1.8504 | 1.9 | 2.0739 | 2 | 1.7151 | 1.8 | 1.4307 | 1.5 | 1.75214 | 1.7 | 1.75214 | 1.7 |
| A6 | 0.9904 | 1 | 0.9809 | 1 | 0.9456 | 1 | 0.8461 | 0.8 | 0.93785 | 0.9 | 0.93785 | 0.9 |
| A7 | 0.5107 | 0.5 | 0.5055 | 0.5 | 0.5812 | 0.5 | 0.5675 | 0.6 | 0.45127 | 0.5 | 0.45127 | 0.5 |
| A8 | 3.4 | 3.4 | 3.4 | 3.4 | 3.3885 | 3.4 | 3.3693 | 3.4 | 3.4493 | 3.5 | 3.4493 | 3.5 |
| Weight (lb) | 484.3278 | 487.07 | 484.3331 | 488.57 | 484.9171 | 485.57 | 489.6039 | 489.49 | 483.8996 | 483.67 | 483.8986 | 483.67 |
| Mean weight (lb) | 491.4399 | | 494.4975 | | 496.2969 | | 516.2156 | | 493.22 | | 492.8937 | |
| Standard Deviation | 13.3153 | | 21.1506 | | 10.2382 | | 16.6668 | | 12.243 | | 9.9664 | |
| Functional Evaluations | 2000 | | 2000 | | 2000 | | 2000 | | | | 2000 | |

4.4 : Benchmark 4 - 52 Bar planar truss

4.4.1    Benchmark 4: Problem description

The 52-bar planar truss is shown in Figure 6. Researchers such as Li et al, Sadollah et al and Do and Lee [34, 36] have studied this problem. The material properties and design constraints of the truss are highlighted in Table 6. Two vertical loads of 10000 lbs each are vertically applied to the truss at nodes 2 and 4. The design variables are categorised into 12 groups as follows: (1) $A_1$–$A_4$; (2) $A_5$–$A_{10}$; (3) $A_{11}$–$A_{13}$,;(4) $A_{14}$–$A_{17}$; (5) $A_{18}$–$A_{23}$; (6) $A_{24}$–$A_{26}$; (7) $A_{27}$–$A_{30}$; (8) $A_{31}$–$A_{36}$; (9)$A_{37}$–$A_{39}$; (10) $A_{40}$–$A_{43}$; (11) $A_{44}$–$A_{49}$, and (12) $A_{50}$–$A_{52}$. The discrete cross-sectional areas are chosen according to the AISC codes presented in Table 14.

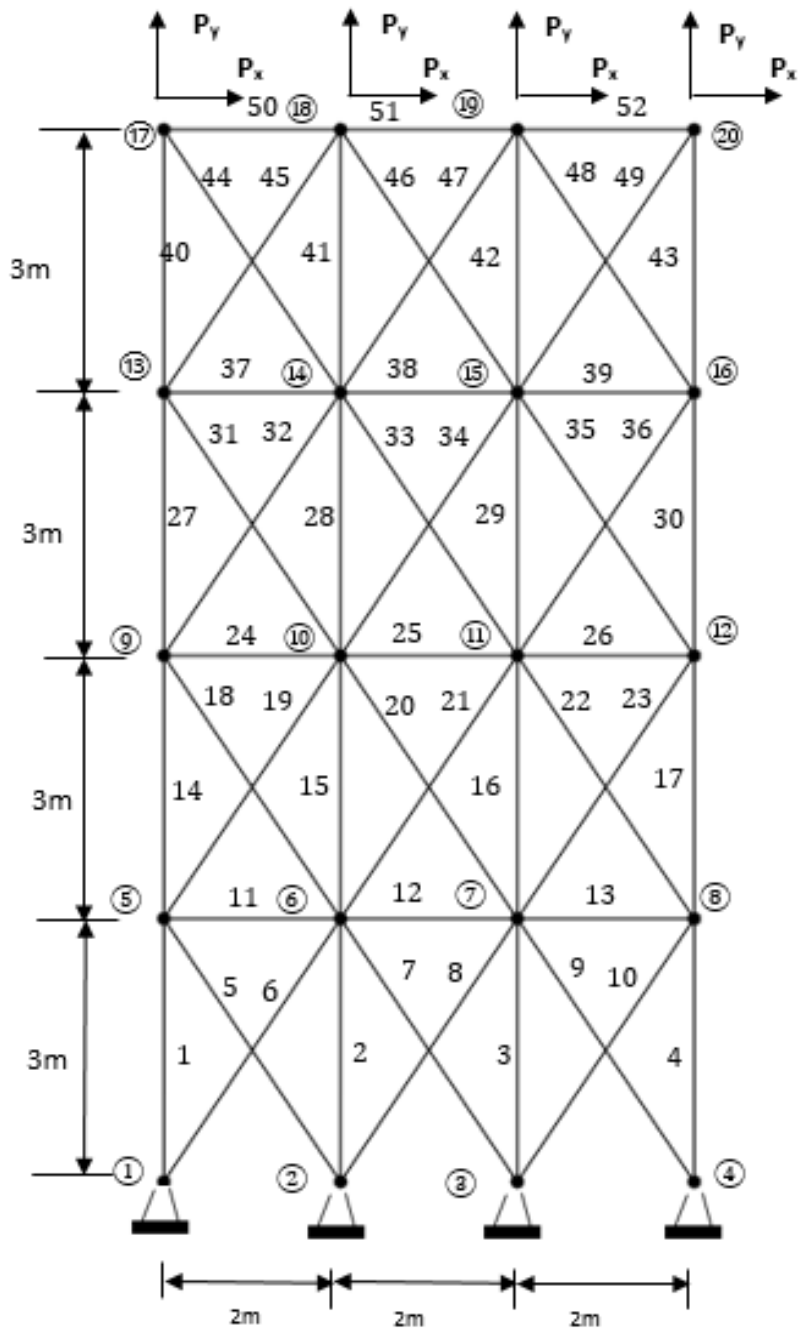

*Figure 6 Benchmark Case 4: The 52-Bar planar truss structure*

*Table 17 The available cross-sectional areas of the AISC code*

| No | in$^2$ | mm$^2$ | No | in$^2$ | mm$^2$ | No | in$^2$ | mm$^2$ | No | in$^2$ | mm$^2$ |
|----|--------|--------|----|--------|--------|----|--------|--------|----|--------|--------|
| 1  | 0.111  | 71.613  | 17 | 1.563 | 1008.385 | 33 | 3.840  | 2477.414 | 49 | 11.500 | 7419.340  |
| 2  | 0.141  | 90.968  | 18 | 1.620 | 1045.159 | 34 | 3.870  | 2496.769 | 50 | 13.500 | 8709.660  |
| 3  | 0.196  | 126.451 | 19 | 1.800 | 1161.288 | 35 | 3.880  | 2503.221 | 51 | 13.900 | 8967.724  |
| 4  | 0.250  | 161.290 | 20 | 1.990 | 1283.868 | 36 | 4.180  | 2696.769 | 52 | 14.200 | 9161.272  |
| 5  | 0.307  | 198.064 | 21 | 2.130 | 1374.191 | 37 | 4.220  | 2722.575 | 53 | 15.500 | 9999.980  |
| 6  | 0.391  | 252.258 | 22 | 2.380 | 1535.481 | 38 | 4.490  | 2896.768 | 54 | 16.000 | 10322.560 |
| 7  | 0.442  | 285.161 | 23 | 2.620 | 1690.319 | 39 | 4.590  | 2961.284 | 55 | 16.900 | 10903.204 |
| 8  | 0.563  | 363.225 | 24 | 2.630 | 1696.771 | 40 | 4.800  | 3096.768 | 56 | 18.800 | 12128.008 |
| 9  | 0.602  | 388.386 | 25 | 2.880 | 1858.061 | 41 | 4.970  | 3206.445 | 57 | 19.900 | 12838.684 |
| 10 | 0.766  | 494.193 | 26 | 2.930 | 1890.319 | 42 | 5.120  | 3303.219 | 58 | 22.000 | 14193.520 |
| 11 | 0.785  | 506.451 | 27 | 3.090 | 1993.544 | 43 | 5.740  | 3703.218 | 59 | 22.900 | 14774.154 |
| 12 | 0.994  | 641.289 | 28 | 3.130 | 2019.351 | 44 | 7.220  | 4658.055 | 60 | 24.500 | 15806.420 |
| 13 | 1.000  | 645.160 | 29 | 3.380 | 2180.641 | 45 | 7.970  | 5141.925 | 61 | 26.500 | 17096.740 |
| 14 | 1.228  | 792.456 | 30 | 3.470 | 2238.705 | 46 | 8.530  | 5503.215 | 62 | 28.000 | 18064.480 |
| 15 | 1.266  | 816.773 | 31 | 3.550 | 2290.318 | 47 | 9.300  | 5999.988 | 63 | 30.000 | 19354.800 |
| 16 | 1.457  | 939.998 | 32 | 3.630 | 2341.931 | 48 | 10.850 | 6999.986 | 64 | 33.500 | 21612.860 |

### 4.4.2    Benchmark 4: Results and discussion

Tables 15 and 16 present the optimal solutions of the 52-bar truss problem. From observation, it can be seen that the TLBO algorithm produces the best truss design with a corresponding weight of 1816.4 kg compared to the s-FOA which obtains a weight of 1819.2kg, DE algorithm with 1820.9 kg, PSO with 1822.5 kg, cFOA with 1839.7 kg and GA with 4207.3 kg. However, from the statistical analysis, the DE algorithm attains the best mean weight of 1830.5 kg and standard deviation of 7.40 kg compared to a mean weight and standard deviation weight of 1834.3 kg in TLBO and 1840.5 kg and 16.73 kg in s-FOA. Hence DE is the most robust and consistent algorithm for solving the 52-bar truss problem of all algorithms implemented.

Nevertheless, from a discrete perspective, the restricted number of structural analysis of 2000 significantly affects the performance of the algorithms (DE, TLBO, PSO, GA, cFOA and s-FOA) investigated in this study as compared that reported in other literatures; 150,000 analyses in HPSO [34], 300,000 analyses in SOS [36] and mSOS [36]. The s-FOA, PSO and TLBO algorithm record the lightest discrete weight of 1912.6 kg which is 1.01% heavier than that obtained by mSOS [36] of 1899.7 kg. The DE algorithm, cFOA and GA produces a weight of 1914.1 kg, 1927.828 kg and 4063.5 kg respectively, 1.01%, 1.02% and 2.14% worse than mSOS [36].

*Table 18 Benchmark Case 4: Performance rankings of the algorithms in solving the 52 bar truss problem*

|  | PERFORMANCE RANKINGS | | | | | |
|---|---|---|---|---|---|---|
|  | DE | PSO | GA | TLBO | cFOA | s-FOA |
| Best Value | 3 | 4 | 6 | 1 | 5 | 2 |
| Mean Value | 1 | 5 | 6 | 2 | 4 | 3 |
| Number of Functional Evaluations per 10,000 | 1.4 | 1.4 | 1.4 | 1.4 | 1.4 | 1.4 |
| Standard deviation | 1 | 5 | 6 | 2 | 4 | 3 |
| Number of tuning parameters changed from one benchmark problem to another | 5 | 3 | 3 | 3 | 3 | 3 |
| Total Score | 11.4 | 16.4 | 19.4 | 9.4 | 17.4 | 12.4 |

*Table 19 Optimal results of the 52-bar planar truss structure optimisation obtained from different algorithms from previous study*

| 52 Bar | | | | | |
|---|---|---|---|---|---|
| | PREVIOUS STUDY | | | | |
| Area group | HPSO [34] | MBA [33] | DE [36] | SOS [36] | mSOS [36] |
| mm$^2$ | Discrete optimum | Discrete optimum | Discrete optimum | Discrete optimum | Discrete optimum |
| A1 | 4658.06 | 4658.06 | 4658.06 | 4658.06 | 4658.06 |
| A2 | 1161.29 | 1161.29 | 1161.29 | 1161.29 | 1161.29 |
| A3 | 363.225 | 494.193 | 494.193 | 494.193 | 494.193 |
| A4 | 3303.22 | 3303.22 | 3303.22 | 3303.22 | 3303.22 |
| A5 | 940 | 940 | 940 | 940 | 940 |
| A6 | 494.193 | 494.193 | 506.451 | 494.193 | 506.451 |
| A7 | 2238.71 | 2238.71 | 2238.71 | 2238.71 | 2238.71 |
| A8 | 1008.39 | 1008.39 | 1008.39 | 1008.39 | 1008.39 |
| A9 | 388.386 | 494.193 | 388.386 | 494.193 | 388.386 |
| A10 | 1283.87 | 1283.87 | 1283.87 | 1283.87 | 1283.87 |
| A11 | 1161.29 | 1161.29 | 1161.29 | 1161.29 | 1161.29 |
| A12 | 792.256 | 494.193 | 506.451 | 494.193 | 506.451 |
| Weight (kg) | 1905.49 | 1902.605 | 1899.654 | 1902.605 | 1899.654 |
| Mean weight | - | - | - | - | - |
| Standard Deviation | - | - | - | - | - |
| Functional Evaluations | 150000 | 5450 | 300,000 | 300,000 | 300,000 |

*Table 20 Optimal results of the 52-bar planar truss structure optimisation obtained from different algorithms in this study*

| 52 Bar | | | | | | | | | | | | |
|---|---|---|---|---|---|---|---|---|---|---|---|---|
| | THIS STUDY | | | | | | | | | | | |
| Area group | TLBO | | PSO | | DE | | GA | | cFOA | | s-FOA | |
| mm² | Continuous Optimum | Discrete optimum | Continuous Optimum | Discrete optimum | Continuous Optimum | Discrete optimum | Continuous Optimum | Discrete optimum | Continuous Optimum | Discrete Optimum | Continuous Optimum | Discrete optimum |
| A1 | 4386.8 | 4658.055 | 4390.6 | 4658.055 | 4396.1 | 4658.055 | 7880.6 | 7419.34 | 4363.48751 | 4658.055 | 4417.6 | 4658.055 |
| A2 | 1129.2 | 1161.288 | 1124.8 | 1161.288 | 1126.7 | 1161.288 | 2008.6 | 1993.544 | 1163.195636 | 1161.288 | 1121.7 | 1161.288 |
| A3 | 318.4 | 252.258 | 282.5 | 252.258 | 293.2 | 285.161 | 837.8 | 792.256 | 369.2652732 | 363.225 | 266.6 | 252.258 |
| A4 | 3376.2 | 3703.218 | 3372.4 | 3703.218 | 3375.6 | 3703.218 | 5822 | 5503.215 | 3445.265624 | 3703.218 | 3393.2 | 3703.218 |
| A5 | 861.9 | 939.998 | 870.3 | 939.998 | 885.5 | 939.998 | 1810.5 | 1696.771 | 878.6127357 | 939.998 | 877.8 | 939.998 |
| A6 | 236 | 252.258 | 271.3 | 252.258 | 252.9 | 252.258 | 1091.6 | 1045.159 | 316.0639755 | 285.161 | 241 | 252.258 |
| A7 | 2295.7 | 2290.318 | 2281.3 | 2290.318 | 2304.2 | 2290.318 | 6212 | 5999.998 | 2292.890711 | 2290.318 | 2312.4 | 2290.318 |
| A8 | 968.3 | 1008.385 | 964.8 | 1008.385 | 964.3 | 1008.385 | 2522.9 | 2503.221 | 981.384552 | 1008.385 | 960.3 | 1008.385 |
| A9 | 282.7 | 285.161 | 294.9 | 285.161 | 268.9 | 285.161 | 3246.8 | 3206.445 | 260.1915129 | 285.161 | 257.3 | 285.161 |
| A10 | 1305 | 1283.868 | 1281.7 | 1283.868 | 1312.5 | 1283.868 | 3926.7 | 3703.218 | 1314.09855 | 1374.191 | 1345.9 | 1283.868 |
| A11 | 1062 | 1161.288 | 1092.8 | 1161.288 | 1068.9 | 1161.288 | 1337.1 | 1696.771 | 1081.592623 | 1045.159 | 1050 | 1161.288 |
| A12 | 451.9 | 494.193 | 500.1 | 494.193 | 433.3 | 494.193 | 8666.8 | 7419.34 | 432.4494231 | 494.193 | 414.4 | 494.193 |
| Weight (kg) | 1816.4 | 1912.524 | 1822.5 | 1912.524 | 1820.9 | 1914.076 | 4207.3 | 4063.478 | 1839.7 | 1927.828 | 1819.2 | 1912.524 |
| Mean weight (kg) | 1834.3 | | 2270.8 | | 1830.5 | | 5295.4 | | 1882.6 | | 1840.5 | |
| Standard Deviation | 15.6369 | | 621.5253 | | 7.4022 | | 875.5519 | | 38.5744 | | 16.7269 | |
| Functional Evaluations | 14000 | | 14000 | | 14000 | | 14000 | | | | 14000 | |

4.5  : Benchmark 5 – 72-Bar planar truss

4.5.1    Benchmark 5: Problem description

The 72-bar planar truss is selected as the fourth problem as shown in Figure 7. The material properties and design constraints of the truss are indicated in Table 6. Table 17 defines the load cases acting on the truss structure. All structural members are categorised into 16 design variables as follows; (1) A1–A4, (2) A5–A12, (3) A13–A16, (4) A17–A18, (5) A19–A22, (6) A23–A30 (7) A31–A34, (8) A35–A36, (9) A37–A40, (10) A41–A48, (11) A49–A52, (12) A53–A54, (13) A55–A58, (14) A59–A66 (15) A67–A70, (16) A71–A72. A discrete set of data as displayed in Table 14 is given for the design.

The truss problem has been previously treated by Li *et al*. and Kaveh and Mahdavi [34, 39].

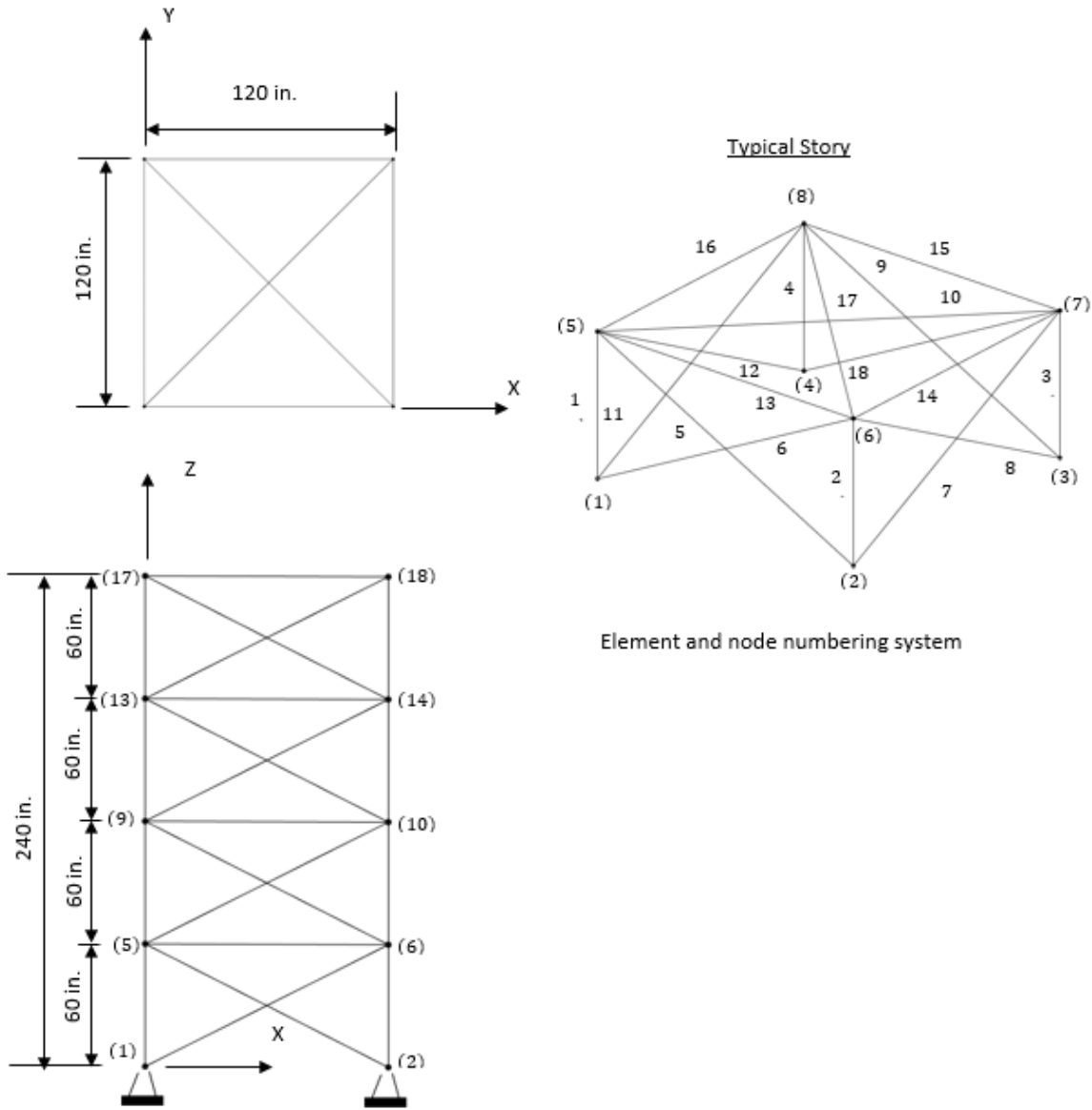

*Figure 7 Benchmark Case 5: The 72-Bar Spatial truss structure*

*Table 21 Benchmark Case 5: Loading Configuration of the 72-bar truss*

| Nodes | Load case 1 | | | Load case 2 | | |
|-------|-------------|-------------|-------------|-------------|-------------|-------------|
| | $P_x$ (kips) | $P_y$ (kips) | $P_z$ (kips) | $P_x$ (kips) | $P_y$ (kips) | $P_z$ (kips) |
| 17 | 5 | 5 | -5 | 0 | 0 | -5 |
| 18 | 0 | 0 | 0 | 0 | 0 | -5 |
| 19 | 0 | 0 | 0 | 0 | 0 | -5 |
| 20 | 0 | 0 | 0 | 0 | 0 | -5 |

### 4.5.2 Benchmark 5: Results and discussion

The optimum 72-bar truss design found by the s-FOA is compared in Tables 18 and 19 to those obtained by other algorithms in this study and other reported literatures. As can be seen, the s-FOA produces yet again the best continuous optimal weight of 403.55lb compared to other algorithms examined. However, the TLBO algorithm attains the best mean weight and standard deviation of 404.12 lb and 0.45 lb respectively compared to the s-FOA algorithm with a mean weight of 406.17 lb and standard deviation of 1.64 lb. The ECBO [39] presents the best discrete truss design with a weight of 389.33 lb. Nevertheless, the number of structural analysis utilised by the ECBO algorithm to achieve the results is 7320 less than that utilised in this study. With a maximum number of functional analysis of 2000, the s-FOA, TLBO, DE, PSO algorithms obtained a discrete weight of 403.22 lb each while the cFOA and GA algorithm produced a weight of 408.51 lb and 594.42 lb respectively.

*Table 22 Benchmark Case 5: Performance rankings of the algorithms in solving the 72-bar truss problem*

| | PERFORMANCE RANKINGS | | | | | |
|---|---|---|---|---|---|---|
| | DE | PSO | GA | TLBO | cFOA | s-FOA |
| Best Value | 4 | 2 | 6 | 3 | 5 | 1 |
| Mean Value | 5 | 4 | 6 | 1 | 3 | 2 |
| Number of Functional Evaluations per 10,000 | 1.05 | 1.05 | 1.05 | 1.05 | 1.05 | 0.96 |
| Standard deviation | 4 | 5 | 6 | 1 | 3 | 2 |
| Number of tuning parameters changed from one benchmark problem to another | 5 | 3 | 3 | 3 | 3 | 3 |
| Total Score | 17.05 | 13.05 | 19.05 | 9.05 | 15.05 | 8.96 |

*Table 23 Optimal results of the 72-bar space truss structure optimisation obtained from different algorithms*

| 72 Bar | | | | | | | | | | | | | | | |
|---|---|---|---|---|---|---|---|---|---|---|---|---|---|---|---|
| **Area group** | PREVIOUS STUDY | | | | **TLBO** | | **PSO** | | THIS STUDY | | | | | | | |
| | SGA [26] | HPSO [22] | CBO [28] | ECBO [28] | | | | | **DE** | | **GA** | | **cFOA** | | **s-FOA** | |
| in$^2$ | Discrete optimum | Discrete optimum | Discrete optimum | Discrete optimum | Continuous Optimum | Discrete optimum | Continuous Optimum | Discrete optimum | Continuous Optimum | Discrete optimum | Continuous Optimum | Discrete optimum | Continuous Optimum | Discrete optimum | Continuous Optimum | Discrete optimum |
| A1 | 0.196 | 4.97 | 1.62 | 1.99 | 1.9691 | 1.8 | 1.9119 | 1.8 | 2.1191 | 1.8 | 2.5972 | 2.38 | 1.9495041 | 1.99 | 1.9877 | 1.8 |
| A2 | 0.602 | 1.228 | 0.563 | 0.563 | 0.5344 | 0.563 | 0.5354 | 0.563 | 0.4829 | 0.563 | 0.3776 | 0.391 | 0.5275275 | 0.563 | 0.5367 | 0.563 |
| A3 | 0.307 | 0.111 | 0.111 | 0.111 | 0.1111 | 0.111 | 0.111 | 0.111 | 0.111 | 0.111 | 0.3366 | 0.307 | 0.1427016 | 0.141 | 0.1086 | 0.111 |
| A4 | 0.766 | 0.111 | 0.111 | 0.111 | 0.1118 | 0.111 | 0.111 | 0.111 | 0.1282 | 0.111 | 0.1524 | 0.141 | 0.1601286 | 0.141 | 0.111 | 0.111 |
| A5 | 0.391 | 2.88 | 1.457 | 1.288 | 1.3882 | 1.457 | 1.3905 | 1.457 | 1.5096 | 1.457 | 1.3759 | 1.266 | 1.3779207 | 1.457 | 1.3584 | 1.457 |
| A6 | 0.391 | 1.457 | 0.442 | 0.442 | 0.5921 | 0.602 | 0.5936 | 0.602 | 0.5753 | 0.602 | 1.116 | 1 | 0.5971467 | 0.563 | 0.6117 | 0.602 |
| A7 | 0.141 | 0.141 | 0.111 | 0.111 | 0.1112 | 0.111 | 0.111 | 0.111 | 0.111 | 0.111 | 0.3481 | 0.307 | 0.1151736 | 0.111 | 0.0882 | 0.111 |
| A8 | 0.111 | 0.111 | 0.111 | 0.111 | 0.1129 | 0.111 | 0.111 | 0.111 | 0.1116 | 0.111 | 0.7811 | 0.766 | 0.1137861 | 0.111 | 0.079 | 0.111 |
| A9 | 0.8 | 1.563 | 0.602 | 0.563 | 0.4895 | 0.442 | 0.5408 | 0.442 | 0.4232 | 0.442 | 2.6031 | 2.38 | 0.5016645 | 0.442 | 0.4876 | 0.442 |
| A10 | 0.602 | 1.228 | 0.563 | 0.563 | 0.5577 | 0.563 | 0.5564 | 0.563 | 0.6054 | 0.563 | 0.4919 | 0.563 | 0.5509707 | 0.563 | 0.5575 | 0.563 |
| A11 | 0.141 | 0.111 | 0.111 | 0.111 | 0.111 | 0.111 | 0.111 | 0.111 | 0.1145 | 0.111 | 0.1121 | 0.111 | 0.1581861 | 0.141 | 0.1106 | 0.111 |
| A12 | 0.307 | 0.196 | 0.111 | 0.111 | 0.111 | 0.111 | 0.111 | 0.111 | 0.128 | 0.111 | 0.1184 | 0.111 | 0.1618047 | 0.141 | 0.1103 | 0.111 |
| A13 | 1.563 | 0.391 | 0.196 | 0.196 | 0.1552 | 0.141 | 0.1543 | 0.141 | 0.1674 | 0.141 | 3.2152 | 3.13 | 0.1502829 | 0.141 | 0.1563 | 0.141 |
| A14 | 0.766 | 1.457 | 0.602 | 0.563 | 0.5713 | 0.563 | 0.5699 | 0.563 | 0.5973 | 0.563 | 0.5636 | 0.563 | 0.5735703 | 0.563 | 0.5803 | 0.563 |
| A15 | 0.141 | 0.766 | 0.391 | 0.391 | 0.408 | 0.391 | 0.4152 | 0.391 | 0.3945 | 0.391 | 0.2581 | 0.307 | 0.3995112 | 0.391 | 0.4351 | 0.391 |
| A16 | 0.111 | 1.563 | 0.563 | 0.563 | 0.5604 | 0.563 | 0.5545 | 0.563 | 0.5077 | 0.563 | 0.6043 | 0.602 | 0.5398152 | 0.563 | 0.4881 | 0.563 |
| Weight (lb) | 427.203 | 393.09 | 391.07 | 389.33 | 404.1224 | 403.22 | 404.0382 | 403.22 | 409.1085 | 403.22 | 615.4444 | 594.42 | 409.1748 | 408.51 | 403.5532 | 403.22 |
| Mean weight (lb) | - | - | - | - | 404.6884 | | 416.7609 | | 420.1844 | | 1268.7 | | 410.6061 | | 406.1689 | |
| Standard Deviation | - | - | - | - | 0.4552 | | 24.263 | | 8.9041 | | 461.4097 | | 1.6644 | | 1.6442 | |
| Functional Evaluations | 60000 | 50000 | 4500 | 3180 | 10500 | | 10500 | | 10500 | | 10500 | | | | 9600 | |

4.6 : Benchmark 6 – 200-Bar planar truss

4.6.1    Benchmark 6: Problem description

Figure 8 illustrates the 200-bar planar truss which has been examined by Lee and Geem [40]. The material properties and design constraints on this example are displayed in Table 6. The members are lined into 29 design categories as shown in Table 20. The minimum possible cross-sectional area is 0.1 in². The truss structure is subjected to three loading conditions as follows: (1) 1 kip in positive x-axis at nodes of 1, 6, 15, 20, 29, 43,48 57, 62 and 71; (2) 10 kips in negative y-axis at nodes of 1, 2, 3, 4,5, 6, 8, 10, 12, 14, 15, 16, 17, 18, 19, 20, 22, 24, 26, 28, 29, 30, 31, 32,33, 34, 36, 38, 40, 42, 43, 44, 45, 46, 47, 48, 50, 52, 54, 56, 58, 59, 60,61, 62, 64, 66, 68, 70, 71, 72, 73, 74 and 75, and (3) both case 1 and case 2.

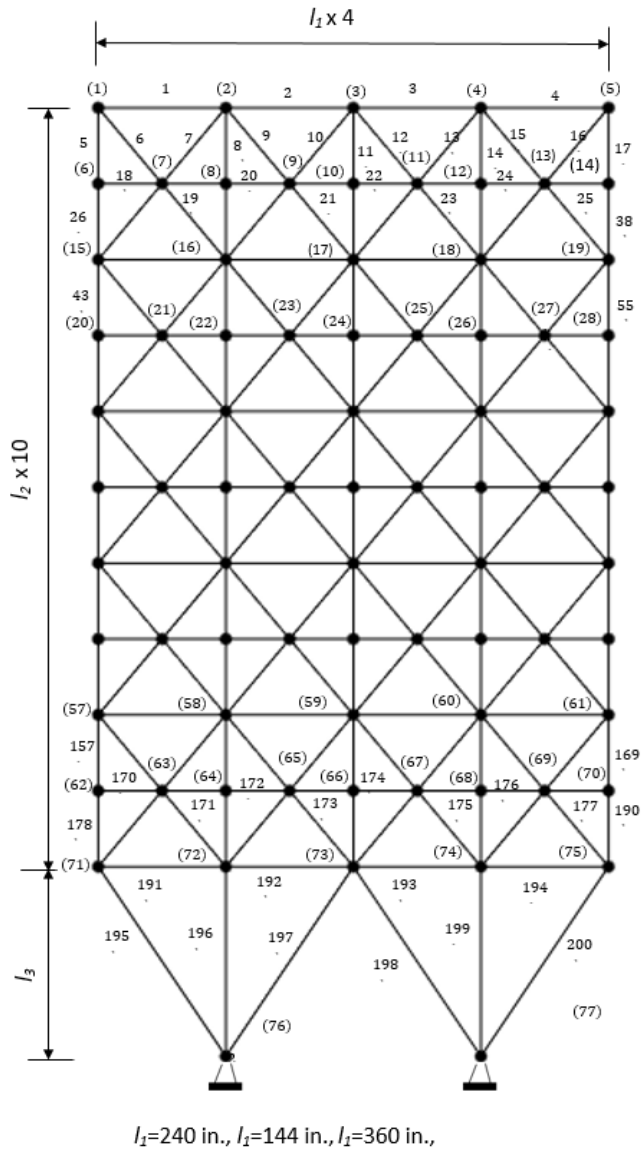

$l_1$=240 in., $l_2$=144 in., $l_3$=360 in.,

*Figure 8 Benchmark Case 6: The 200-Bar planar truss structure*

*Table 24 Benchmark Case 6: Design variable classification of the 200-bar truss*

| Element Group Number | Members in the group |
|---|---|
| 1 | 1,2,3,4 |
| 2 | 5,8,11,14,17 |
| 3 | 19,20,21,22,23,24 |
| 4 | 18,25, 56, 63, 94, 101, 132, 139, 170, 177 |
| 5 | 26, 29, 32, 35, 38 |
| 6 | 6,7,9,10, 12, 13, 15, 16, 27, 28, 30, 31, 33, 34, 36, 37 |
| 7 | 39, 40, 41, 42 |
| 8 | 43, 46, 49, 52, 55 |
| 9 | 57, 58, 59, 60, 61, 62 |
| 10 | 64, 67, 70, 73, 76 |
| 11 | 44, 45, 47, 48, 50, 51, 53, 54, 65, 66, 68, 69, 71, 72, 74, 75 |
| 12 | 77, 78, 79, 80 |
| 13 | 81, 84, 87, 90, 93 |
| 14 | 95, 96, 97, 98, 99, 100 |
| 15 | 102, 105, 108, 111, 114 |
| 16 | 82, 83, 85, 86, 88, 89, 91, 92, 103, 104, 106, 107, 109, 110, 112, 113 |
| 17 | 115, 116, 117, 118 |
| 18 | 119, 122, 125, 128, 131 |
| 19 | 133, 134, 135, 136, 137, 138 |
| 20 | 140, 143, 146, 149, 152 |
| 21 | 120, 121, 123, 124, 126, 127, 129, 130, 141, 142, 144, 145, 147, 148, 150, 151 |
| 22 | 153, 154, 155, 156 |
| 23 | 157, 160, 163, 166, 169 |
| 24 | 171, 172, 173, 174, 175, 176 |
| 25 | 178, 181, 184, 187 190 |
| 26 | 158, 159, 161, 162, 164, 165, 167, 168, 179, 180, 182, 183, 185, 186, 188, 189 |
| 27 | 191, 192, 193, 194 |
| 28 | 195, 197, 198, 200 |
| 29 | 196, 199 |

### 4.6.2    Benchmark 6: Results and discussion

Tables 21 and 22 display the results obtained from the considered algorithms with other optimisation techniques. The results show that the s-FOA algorithm proposes a truss design with a weight of 25,505lb. The HS based algorithm [40] ranks first amongst the algorithms considered in the study with a weight of 25,447.1lb. However, the number of functional evaluations utilised by the HS based algorithm [40] is 3000 more than the other algorithms considered. Nonetheless, the s-FOA algorithm is placed third in the comparison. The truss weights obtained from the TLBO, DE, cFOA GA and PSO algorithms are 25497 lb, 26189 lb, 26746 lb, 31264 lb, and 25962 lb respectively. The result also indicates that the s-FOA compared to other algorithms is placed third in terms of the best mean weights and standard deviation with values of 26891lb and 537lb respectively. The HS based algorithm [40] is exempted from this analysis as there is no information reported on the mean weight and standard deviation to establish a comparison.

*Table 25 Benchmark Case 6: Performance rankings of the algorithms in solving the 200-bar truss problem*

| | PERFORMANCE RANKINGS | | | | | |
|---|---|---|---|---|---|---|
| | DE | PSO | GA | TLBO | cFOA | s-FOA |
| Best Value | 4 | 3 | 6 | 1 | 5 | 2 |
| Mean Value | 2 | 5 | 6 | 1 | 4 | 3 |
| Number of Functional Evaluations per 10,000 | 4.5 | 4.5 | 4.5 | 4.5 | 4.5 | 4.5 |

| Standard deviation | 1 | 5 | 6 | 2 | 4 | 3 |
|---|---|---|---|---|---|---|
| Number of tuning parameters changed from one benchmark problem to another | 5 | 3 | 3 | 3 | | 3 |
| Total Score | 16.5 | 17.5 | 21.5 | 11.5 | 17.5 | 12.5 |

*Table 26 Optimal results of the 200-bar planar truss structure optimisation obtained from different algorithms*

| | 200 Bar | | | | | | |
|---|---|---|---|---|---|---|---|
| | PREVIOUS STUDY | | THIS STUDY | | | | |
| | Lee and Geem [40] | TLBO | PSO | DE | GA | cFOA | s-FOA |
| **Area group** in² | Continuous Optimum | Continuous Optimum | Continuous Optimum | Continuous Optimum | Continuous Optimum | Continuous Optimum | Continuous Optimum |
| A1 | 0.1253 | 0.1002 | 0.1 | 0.1464 | 0.1976 | 0.84276 | 0.1008 |
| A2 | 1.0157 | 0.9446 | 0.9763 | 1.1397 | 1.0795 | 1.42356 | 0.9946 |
| A3 | 0.1069 | 0.35 | 0.1 | 0.2703 | 0.1196 | 0.1783 | 0.1035 |
| A4 | 0.1096 | 0.1124 | 0.1 | 0.1481 | 0.1992 | 0.97888 | 0.1184 |
| A5 | 1.9369 | 1.9457 | 1.9781 | 2.2242 | 4.5599 | 2.00472 | 1.9909 |
| A6 | 0.2686 | 0.2879 | 0.145 | 0.2667 | 0.1386 | 0.30422 | 0.1651 |
| A7 | 0.1042 | 0.1449 | 0.4787 | 0.187 | 0.1442 | 0.44632 | 0.3790 |
| A8 | 2.9731 | 3.1724 | 3.0254 | 3.2557 | 5.8761 | 3.20697 | 3.0480 |
| A9 | 0.1309 | 0.1011 | 0.114 | 0.1 | 0.2381 | 0.20619 | 0.1100 |
| A10 | 4.1831 | 4.1576 | 4.049 | 4.2936 | 4.5466 | 4.13978 | 4.0586 |
| A11 | 0.3967 | 0.3097 | 0.3991 | 0.3697 | 0.2686 | 0.41515 | 0.5880 |
| A12 | 0.4416 | 0.1824 | 0.4129 | 0.3315 | 0.2601 | 0.17381 | 0.2053 |
| A13 | 5.1873 | 5.3714 | 5.2903 | 5.5847 | 4.9692 | 5.29643 | 5.3936 |
| A14 | 0.1912 | 0.1417 | 0.4163 | 0.1 | 1.1107 | 0.83119 | 0.2111 |
| A15 | 6.241 | 6.4221 | 6.2992 | 6.481 | 5.7642 | 6.25224 | 6.3832 |
| A16 | 0.6994 | 0.4274 | 0.7178 | 0.6262 | 0.9245 | 0.79046 | 0.4990 |
| A17 | 0.1158 | 0.548 | 0.3752 | 0.4016 | 0.5035 | 0.13696 | 0.3609 |
| A18 | 7.7643 | 7.7648 | 8.0483 | 8.0576 | 7.1131 | 8.64241 | 7.8628 |
| A19 | 0.1 | 0.1099 | 0.6505 | 0.1297 | 0.9186 | 0.12575 | 0.1003 |
| A20 | 8.8279 | 8.7661 | 9.0707 | 8.9544 | 7.9984 | 8.98813 | 8.8603 |
| A21 | 0.6986 | 0.7582 | 1.1272 | 0.7092 | 0.6227 | 0.6403 | 0.6673 |
| A22 | 1.5563 | 0.502 | 0.2648 | 1.2072 | 5.5998 | 0.20019 | 0.4757 |
| A23 | 10.9806 | 10.6533 | 11.7034 | 11.3708 | 8.6013 | 10.60984 | 10.6563 |

| A24 | 0.1317 | 0.6135 | 0.9308 | 0.3314 | 1.8583 | 0.33525 | 0.2405 |
|---|---|---|---|---|---|---|---|
| A25 | 12.1492 | 11.6602 | 12.7045 | 12.4438 | 9.6234 | 11.65307 | 11.6516 |
| A26 | 1.6373 | 1.3 | 1.6645 | 1.6162 | 2.3529 | 0.86093 | 1.0313 |
| A27 | 5.0032 | 6.4437 | 4.1803 | 5.0998 | 10.7055 | 6.94793 | 6.9146 |
| A28 | 9.3545 | 10.5826 | 9.0751 | 10.0299 | 14.4036 | 11.96894 | 10.8883 |
| A29 | 15.091 | 13.9279 | 15.4446 | 14.6571 | 12.9477 | 13.27711 | 13.5963 |
| Weight (lb) | 25447.1 | 25497 | 25962 | 26189 | 31264 | 26746 | 25505 |
| Mean weight (lb) | - | 26173 | 28655 | 26735 | 41393 | 27727 | 26891 |
| Standard Deviation | - | 516.6024 | 18526 | 267.126 | 52944 | 766 | 537 |
| Functional Evaluations | 48000 | 45000 | 45000 | 45000 | 45000 | | 45000 |

## 5. Conclusion

In this study, the Spontaneous Fruit-fly Optimisation Algorithm (s-FOA) is presented to design optimal trusses, subject to stress and displacement bounds. An improved vision phase is proposed as an improvement to the FOA, in order to improve the exploitative capabilities of the algorithm in the search space. The improved vision phase aims to improve control of exploration and exploitation. Standard tuning parameters for use in any structural design are also presented. This eliminates the problems associated with the selection of the right parameters for the s-FOA algorithm. The effectiveness of the s-FOA algorithm is tested on six benchmark truss problems of size optimisation. The design variable is selected as the cross-sectional area of the truss members.

This study compares the performance of the s-FOA algorithm with the cFOA, TLBO, DE, PSO and the GA algorithms as well as other algorithms reported in several literatures such as mSOS, SOS, MBA, HSPO, HS, PSOPC, aeDE and SGA. It is observed that in most of the truss problems the s-FOA obtains the best continuous results in terms of the lowest weight, mean weight and standard deviation compared to other algorithms, despite a restricted number of structural analysis. The s-FOA is also seen to be competitive in discovering discrete optimal designs to the truss problems notwithstanding restricted computing resources. Therefore, s-FOA can be employed by industries with limited computation power in designing optimal truss structures.

## References

1. Maskulus, M. and Schafhirt, S. (2014) "Design Optimization Of Wind Turbine Support Structures-A Review". Journal Of Ocean And Wind Energy 1, 12-22
2. C, L. and Y, O. (2016) "Structural Analysis of Lattice Steel Transmission Towers: A Review". Journal Of Steel Structures & Construction 2 (1)
3. Do, D. and Lee, J. (2017) "A Modified Symbiotic Organisms Search (Msos) Algorithm for Optimization of Pin-Jointed Structures". Applied Soft Computing61, 683-699
4. Felkner, J., Chatzi, E. and Kotnik, T. (2015) "Interactive Truss Design Using Particle Swarm Optimization and NURBS Curves". Journal Of Building Engineering4, 60-74
5. Onyekpe, Uche & Kanarachos, Stratis. (2018). Truss Optimisation: A comparison of the improved contrast-based fruit fly optimisation algorithm (cFOA) with popular swarm optimisation methods. 10.13140/RG.2.2.15820.72326.
6. Kaveh, A. and Mahdavi, V. (2018). Multi-objective colliding bodies optimization algorithm for design of trusses. Journal of Computational Design and Engineering.
7. Assimi, H., Jamali, A. and Nariman-zadeh, N. (2017). Sizing and topology optimization of truss structures using genetic programming. Swarm and Evolutionary Computation, 37, pp.90-103.
8. Cheng, M., Prayogo, D., Wu, Y. and Lukito, M. (2016). A Hybrid Harmony Search algorithm for discrete sizing optimization of truss structure. Automation in Construction, 69, pp.21-33.

9. Tejani, G., Savsani, V., Patel, V. and Savsani, P. (2018). Size, shape, and topology optimization of planar and space trusses using mutation-based improved metaheuristics. Journal of Computational Design and Engineering, 5(2), pp.198-214.

10. Tejani, G., Savsani, V., Patel, V. and Mirjalili, S. (2018). Truss optimization with natural frequency bounds using improved symbiotic organisms search. Knowledge-Based Systems, 143, pp.162-178.

11. Kaveh, A. and Zolghadr, A. (2014) "Democratic PSO for Truss Layout and Size Optimization with Frequency Constraints". Computers & Structures 130, 10-21

12. Farshchin, M., Camp, C. and Maniat, M. (2016) "Multi-Class Teaching–Learning-Based Optimization For Truss Design With Frequency Constraints". Engineering Structures 106, 355-369

13. Rajan, S. (1995) "Sizing, Shape, And Topology Design Optimization Of Trusses Using Genetic Algorithm". Journal Of Structural Engineering 121 (10), 1480-1487

14. Pan, W. (2012). A new Fruit Fly Optimization Algorithm: Taking the financial distress model as an example. Knowledge-Based Systems, 26, pp.69-74.

15. Iscan, H. and Gunduz, M. (2017). An application of fruit fly optimization algorithm for traveling salesman problem. Procedia Computer Science, 111, pp.58-63.

16. Lu, M., Zhou, Y., Luo, Q. and Huang, K. (2018). An Adaptive Fruit Fly Optimization Algorithm Based on Velocity Variable.

17. Kanarachos, S., Griffin, J. and Fitzpatrick, M. (2017). Efficient truss optimization using the contrast-based fruit fly optimization algorithm. Computers & Structures, 182, pp.137-148.

18. Kanarachos, S., Dizqah, A., Chrysakis, G. and Fitzpatrick, M. (2018). Optimal design of a quadratic parameter varying vehicle suspension system using contrast-based Fruit Fly Optimisation. *Applied Soft Computing*, 62, pp.463-477.

19. Kanarachos, S., Savitski, D., Lagaros, N. & Fitzpatrick, M.E. 2017, "Automotive magnetorheological dampers: modelling and parameter identification using contrast-based fruit fly optimisation", Soft Computing, , pp. 1-19.

20. Mitić, M., Vuković, N., Petrović, M. and Miljković, Z. (2015). Chaotic fruit fly optimization algorithm. Knowledge-Based Systems, 89, pp.446-458.

21. Yuan, X., Dai, X., Zhao, J. and He, Q. (2014). On a novel multi-swarm fruit fly optimization algorithm and its application. Applied Mathematics and Computation, 233, pp.260-271.

22. Liu, C., Huang, G., Chai, Q. and Zhang, R. (2014). A modified fruit fly optimization algorithm with better balance between exploration and exploitation. 2014 IEEE International Conference on Progress in Informatics and Computing.

23. Van Breugel, F. and Dickinson, M. (2014). Plume-Tracking Behavior of Flying Drosophila Emerges from a Set of Distinct Sensory-Motor Reflexes. Current Biology, 24(3), pp.274-286.

24. Magdeburg (2018). Memories of fruit fly larvae are more complex than thought — Bernstein Netzwerk Computational Neuroscience. [online] Bernstein-network.de. Available at: http://www.bernstein-network.de/en/news/Forschungsergebnisse-en/memories-of-fruit-fly-larvae-are-more-complex-than-thought [Accessed 21 Apr. 2018].

25. Heris, S. (2015). Implementation of Differential Evolution (DE). [online] Yarpiz.com. Available at: http://yarpiz.com/231/ypea107-differential-evolution [Accessed 19 Sep. 2017].

26. Heris, S. (2015). Implementation of Particle Swarm Optimization in MATLAB. [online] Yarpiz.com. Available at: http://yarpiz.com/50/ypea102-particle-swarm-optimization [Accessed 19 Sep. 2017].

27. Heris, S. (2017). Implementation of TLBO in MATLAB. [online] Yarpiz.com. Available at: http://yarpiz.com/83/ypea111-teaching-learning-based-optimization [Accessed 19 Sep. 2017].

28. Matlab Optimisation Toolbox. Natick, Massachusetts, United States: The MathWorks, Inc., 2017. Print.

29. Magnus Erik Hvass Pedersen (2010). Good Parameters for Differential Evolution. Technical Report no. HL1002. [online] Hvass Laboratories, pp.1-10. Available at: https://pdfs.semanticscholar.org/48aa/36e1496c56904f9f6dfc15323e0c45e34a4c.pdf [Accessed 19 Sep. 2017].

30. Shi, Y. and Eberhart, R. (1998). Parameter selection in particle swarm optimization. In: International Conference on Evolutionary Programming, Evolutionary Programming VII. [online] Springer Link, pp.591-600. Available at: https://link.springer.com/chapter/10.1007/BFb0040810 [Accessed 19 Sep. 2017].

31. Alajmi, A. and Wright, J. (2014). Selecting the most efficient genetic algorithm sets in solving unconstrained building optimization problem. International Journal of Sustainable Built Environment, [online] 3(1), pp.18-26. Available at: http://www.sciencedirect.com/science/article/pii/S2212609014000399 [Accessed 19 Sep. 2017].

32. Farshchin, M., Camp, C. and Maniat, M. (2016). Optimal design of truss structures for size and shape with frequency constraints using a collaborative optimization strategy. Expert Systems with Applications, 66, pp.203-218.

33. Sadollah, Ali et al. "Mine Blast Algorithm for Optimization of Truss Structures with Discrete Variables." Computers & Structures 102-103 (2012): 49-63. Web. 14 Mar. 2018.

34. Li, L.J., Z.B. Huang, and F. Liu. "A Heuristic Particle Swarm Optimization Method for Truss Structures with Discrete Variables." Computers & Structures 87.7-8 (2009): 435-443. Web. 14 Mar. 2018.

35. Ho-Huu, V. et al. "An Adaptive Elitist Differential Evolution for Optimization of Truss Structures With Discrete Design Variables." Computers & Structures 165 (2016): 59-75. Web. 14 Mar. 2018.

36. Do, D. and Lee, J. (2017). A modified symbiotic organisms search (mSOS) algorithm for optimization of pin-jointed structures. Applied Soft Computing, [online] 61, pp.687-697. Available at: http://www.sciencedirect.com/science/article/pii/S1568494617304842 [Accessed 19 Sep. 2017].

37. Wu, Shyue-Jian, and Pei-Tse Chow. "Steady-State Genetic Algorithms for Discrete Optimization of Trusses." Computers & Structures 56.6 (1995): 979-991. Web. 14 Mar. 2018.

38. Lee, Kang Seok et al. "The Harmony Search Heuristic Algorithm for Discrete Structural Optimization." Engineering Optimization 37.7 (2005): 663-684. Web. 14 Mar. 2018.

39. Lee, Kang Seok et al. "The Harmony Search Heuristic Algorithm for Discrete Structural Optimization." Engineering Optimization 37.7 (2005): 663-684. Web. 14 Mar. 2018.
40. Lee, Kang Seok, and Zong Woo Geem. "A New Structural Optimization Method Based On The Harmony Search Algorithm." Computers & Structures 82.9-10 (2004): 781-798. Web. 14 Mar. 2018.